\documentclass[aos,MSNbibl,seceqn,rotating,dvips]{arximspdf}
\usepackage{mathrsfs,dcolumn,multirow}

\doi{10.1214/13-AOS1183} 
\volume{42}
\issue{1}
\pubyear{2014}
\firstpage{142}
\lastpage{170}

\makeatletter
\newcolumntype{d}[1]{D{.}{.}{#1}}

\newcommand{\rright}{\right}
\newcommand{\lleft}{\left}
\newtheorem{theorem}{Theorem}[section]
\newtheorem{corollary}{Corollary}[section]
\newproclaim{remark}{Remark}[section]
\newproclaim{ass}{Assumption}
\newtheorem{proposition}{Proposition}[section]
\let\widebar\overline
\makeatother

\begin{document}
\begin{frontmatter}

\title{Theory and methods of panel data models with~interactive effects}
\runtitle{Panel data models with interactive effects}

\begin{aug}
\author[A]{\fnms{Jushan} \snm{Bai}\corref{}\ead[label=e1]{jushan.bai@columbia.edu}\thanksref{t1}}
\and
\author[B]{\fnms{Kunpeng} \snm{Li}\ead[label=e3]{likp.07@sem.tsinghua.edu.cn}\ead[label=u1,url]{http://www.foo.com}\thanksref{t2}}
\runauthor{J. Bai and K. Li}
\affiliation{Columbia University and Central University of Economics and Finance,\\
Capital University of Economics and Business and Tsinghua University}
\address[A]{Department of Economics\\
Columbia University\\
420, West 118th Street\\
New York, New York 10027\\
USA\\
\printead{e1}} 
\address[B]{International School of Economics and Management\\
Capital University of Economics and Business\\
Huaxiang Town, Fengtai District\\
Beijing 100070\\
China\\
\printead{e3}}
\end{aug}
\thankstext{t1}{Supported by the NSF (SES-0962410).}
\thankstext{t2}{Supported by NSFC (71201031) and Humanities and Social Sciences of Chinese Ministry of Education (12YJCZH109).}

\received{\smonth{12} \syear{2012}}
\revised{\smonth{10} \syear{2013}}

%
\begin{abstract}
This paper considers the maximum likelihood estimation of panel data
models with interactive effects.
Motivated by applications in economics and other social sciences, a
notable feature of the model is
that the explanatory variables are correlated with the unobserved effects.
The usual within-group estimator is inconsistent. Existing methods for
consistent estimation are
either designed for panel data with short time periods or are less
efficient. The maximum likelihood
estimator has desirable properties and is easy to implement, as
illustrated by the Monte Carlo
simulations. This paper develops the inferential theory for the maximum
likelihood estimator,
including consistency, rate of convergence and the limiting
distributions. We further extend
the model to include time-invariant regressors and common regressors
(cross-section invariant).
The regression coefficients for the time-invariant regressors are
time-varying, and the
coefficients for the common regressors are cross-sectionally varying.
\end{abstract}

%

\end{frontmatter}

\section{Introduction}\label{sec1}
This paper studies the following panel data models with unobservable
interactive effects:
\[
y_{it}=\alpha_i+ x_{it}\beta+
\lambda_i'f_t + e_{it}, \qquad
i=1,\ldots,N, t=1,2,\ldots,T;
\]
where $y_{it}$ is the dependent variable;
$x_{it}=(x_{it1},\ldots,x_{itK})$ is a row vector of explanatory
variables; $\alpha_i$ is an intercept; the term $\lambda_i'f_t+e_{it}$
is unobservable and has a factor structure, $\lambda_i$ is an $r\times
1$ vector of factor loadings, $f_t$ is a vector of factors and $e_{it}$
is the idiosyncratic error. The interactive effects ($\lambda_i'f_t$)
generalize the usual additive individual and time effects; for example,
if $\lambda_i\equiv1$, then $\alpha_i+\lambda_i' f_t = \alpha_i+f_t$.

A key feature of the model is that the regressors $x_{it}$ are allowed
to be correlated with $(\alpha_i, \lambda_i, f_t)$. This situation is
commonly encountered in economics and other social sciences, in which
some of the regressors $x_{it}$ are decision variables that are
influenced by the unobserved individual heterogeneities. The practical
relevance of the model will be further discussed below. The objective
of this paper is to obtain consistent and efficient estimation of
$\beta
$ in the presence of correlations between the regressors and the factor
loadings and factors.

The usual pooled least squares estimator or even the within-group
estimator is inconsistent for $\beta$. One method to obtain a
consistent estimator is to treat $(\alpha_i,\lambda_i,f_t)$ as
parameters and estimate them jointly with $\beta$. The idea is
``controlling through estimating'' (controlling the effects by
estimating them). This is the approach used in \cite{bai2009panel,moon2009likelihood} and \cite{su2012specification}. While there
are some advantages, an undesirable consequence of this approach is the
incidental parameters problem. There are too many parameters being
estimated, and the incidental parameters bias arises; see \cite
{neyman1948consistent}. In \cite{ahn2001gmm,ahn2006panel} and
\cite{holtz1988estimating} the authors consider the generalized method
of moments (GMM) method. The GMM method is based on a nonlinear
transformation known as quasi-differencing that eliminates the factor
errors. Quasi-differencing increases the nonlinearity of the model
especially with more than one factor. The GMM method works well with a
small $T$. When $T$ is large, the number of moment equations will be
large, and the so called many-moment bias arises. In \cite{pesaran2006estimation},
the author considers an alternative method by augmenting
the model with additional regressors $\bar y_t$ and $\bar x_t$, which
are the cross-sectional averages of $y_{it}$ and $x_{it}$. These
averages provide an estimate for $f_t$. The estimator of \cite
{pesaran2006estimation} becomes inconsistent when the factor loadings
in the $y$ equation are correlated with those in the $x$ equation, as
shown in \cite{westerlundurbain2013}. A further approach to
controlling the correlation between the regressors and factor errors is
to use the Mundlak--Chamberlain projection (\cite{mundlak1978on}
and~\cite{chamberlain1984panel}). The latter method projects $\alpha_i$ and
$\lambda_i$ onto the regressors such that $\lambda_i= c_0+c_1
x_{i1}+\cdots+c_T x_{iT}+\eta_i$, where $c_s$ ($s=0,1,\ldots,T)$ are
parameters to be estimated, and $\eta_i$ is the projection residual (a
similar projection is done for~$\alpha_i$). The projection residuals
are uncorrelated with the regressors so that a variety of approaches
can be used to estimate the model. This framework is designed for small
$T$ and is studied by \cite{bai2009likelihood}.

In this paper we consider the pseudo-Gaussian maximum likelihood method
under large $N$ and large $T$. The theory does not depend on normality.
In view of the importance of the MLE in the statistical literature, it
is of both practical and theoretical interest to examine the MLE in
this context. We develop a rigorous theory for the MLE. We show that
there is no incidental parameters bias for $\beta$.

We allow time-invariant regressors such as education, race and gender
in the model. The corresponding regression coefficients are
time-dependent. Similarly, we allow common regressors, which do not
vary across individuals, such as prices and policy variables. The
corresponding regression coefficients are individual-dependent so that
individuals respond differently to policy or price changes. In our
view, this is a sensible way to incorporate time-invariant and common
regressors. For example, wages associated with education and with
gender are more likely to change over time rather than remain constant.
In our analysis, time invariant regressors are treated as the
components of $\lambda_i$ that are observable, and common regressors as
the components of $f_t$ that are observable. This view fits naturally
into the factor framework in which part of the factor loadings and
factors are observable, and the maximum likelihood method imposes the
corresponding loadings and factors at their observed values.

While the theoretical analysis of MLE is demanding, the limiting
distributions of the MLE are simple and have intuitive interpretations.
The computation is also easy and can be implemented by adapting the ECM
(expectation and constrained maximization) of \cite{meng1993maximum}.
In addition, the maximum likelihood method allows restrictions to be
imposed on $\lambda_i$ or on $f_t$ to achieve more efficient
estimation. These restrictions can take the form of known values, being
either zeros, or other fixed values. Part of the rigorous analysis
includes setting up the constrained maximization as a Lagrange
multiplier problem. This approach provides insight into which kinds of
restrictions provide efficiency gain and which kinds do not.

Panel data models with interactive effects have wide applicability in
economics. In macroeconomics, for example, $y_{it}$ can be the output
growth rate for country $i$ in year $t$; $x_{it}$ represents production
inputs, and $f_t$ is a vector of common shocks (technological progress,
financial crises); the common shocks have heterogenous impacts across
countries through the different factor loadings $\lambda_i$; $e_{it}$
represents the country-specific unmeasured growth rates. In
microeconomics, and especially in earnings studies, $y_{it}$ is the
wage rate for individual $i$ for period $t$ (or for cohort $t$),
$x_{it}$ is a vector of observable characteristics such as marital
status and experience; $\lambda_i$ is a vector of unobservable
individual traits such as ability, perseverance, motivation and
dedication; the payoff to these individual traits is not constant over
time, but time varying through $f_t$; and $e_{it}$ is idiosyncratic
variations in the wage rates. In finance, $y_{it}$ is stock $i$'s
return in period $t$, $x_{it}$ is a vector of observable factors, $f_t$
is a vector of unobservable common factors (systematic risks) and
$\lambda_i$ is the exposure to the risks; $e_{it}$ is the idiosyncratic
returns. Factor error structures are also used as a flexible trend
modeling as in \cite{kneip2012new}. Most of panel data analysis assumes
cross-sectional independence; see, for example, \cite
{arellano2003panel,baltagi2005econometric} and \cite
{hsiao2003analysis}. The factor structure is also capable of capturing
the cross-sectional dependence arising from the common shocks $f_t$. Further motivation can be found in
\cite{bai2003inferential,ross1976arbitrage,stock2002forecasting}.\vadjust{\goodbreak}

Throughout the paper, the norm of a vector or matrix is that of
Frobenius, that is, $\|A\|=[\operatorname{tr}(A'A)]^{1/2}$ for matrix
$A$; $\operatorname{diag}(A)$
is a column vector consisting of the diagonal elements of $A$ when $A$
is matrix, but $\operatorname{diag}(A)$ represents a diagonal matrix when $A$ is a
vector. In addition, we use $\dot v_t$ to denote $v_t-\frac{1} T\sum_{t=1}^Tv_t$
for any column vector $v_t$ and $M_{wv}$ to denote $\frac{1} T\sum_{t=1}^T\dot w_t\dot v_t'$ for any vectors $w_t$ and $ v_t$.

The rest of the paper is organized as follows. Section~\ref{sec2} introduces a
common shock model and the maximum likelihood estimation. Consistency,
rate of convergence and the limiting distributions of the MLE are
established. Section~\ref{sec3} shows that if some factors do not affect the $y$
equation but only the $x$ equation, more efficient estimation can be
obtained. Section~\ref{sec4} extends the analysis to time-invariant regressors
and common regressors; the corresponding coefficients are time varying
and cross-section varying, respectively.
Computing algorithm is discussed in Section~\ref{sec5}, and simulations results
are reported in Section~\ref{sec6}. The last section concludes. The theoretical
proofs are provided in the supplementary document \cite{bai}.

\section{A common shock model}\label{sec2}

In the common-shock model, we assume that both $y_{it}$ and $x_{it}$
are impacted by the common shocks $f_t$ so the model takes the form
\begin{eqnarray}
\label{comshkmaineq1} y_{it}&=&\alpha_i+x_{it1}
\beta_1+x_{it2}\beta_2+\cdots +x_{itK}
\beta_K+\lambda_i'f_t+e_{it},
\nonumber\\[-8pt]\\[-8pt]
x_{itk}&=&\mu_{ik}+\gamma_{ik}'f_t+v_{itk}\nonumber
\end{eqnarray}
for $k=1,2,\ldots,K$. In across-country output studies, for example,
output $y_{it}$ and inputs $x_{it}$ (labor and capital) are both
affected by the common shocks.

The parameter of interest is $\beta=(\beta_1,\ldots,\beta_K)'$. We also estimate
$\alpha_i,\lambda_i,\mu_{ik}$ and $\gamma_{ik}$ $(k=1,2,\ldots,K)$. By
treating the latter as parameters, we also allow arbitrary correlations
between $(\alpha_i,\lambda_i)$ and $(\mu_{ik},\gamma_{ik})$. Although
we also treat $f_t$ as fixed parameters, there is no need to estimate
the individual $f_t$, but only the sample covariance of $f_t$. This is
an advantage of the maximum likelihood method, which eliminates the
incidental parameters problem in the time dimension. This kind of the
maximum likelihood method was used for pure factor models in \cite{amemiya1987asymptotic,amemiya88} and~\cite{bai2012statistical}.
By symmetry, we could also estimate individuals $f_t$, but then we only
estimate the sample covariance of the factor loadings. The idea is that
we do not simultaneously estimate the factor loadings and the factors
$f_t$ (which would be the case for the principal components method).
This reduces the number of parameters considerably. If $N$ is much
smaller than $T$ $(N\ll T)$, treating factor loadings as parameters is
preferable since there are fewer parameters.

Because of the correlation between the regressors and regression errors
in the $y$ equation, the $y$ and $x$ equations form a simultaneous
equation system; the MLE jointly estimates the parameters in both
equations. The joint estimation avoids the Mundlak--Chamberlain
projection and thus is applicable for large $N$ and large $T$.

We assume the number of factors $r$ is fixed and known. Determining the
number of factors is discussed in Section~\ref{sec6}, where a modified
information criterion proposed by \cite{bai2002determining} is used.
Let $x_{it}=(x_{it1},x_{it2},\ldots,x_{itK})$, $\gamma
_{ix}=(\gamma
_{i1},\gamma_{i2},\dots,\gamma_{iK})$,
$v_{itx}=(v_{it1},v_{it2},\ldots,v_{itK})'$ and $\mu_{ix}=(\mu_{i1},\mu_{i2},\ldots,\mu_{iK})'$. The
second equation of (\ref{comshkmaineq1}) can be written in matrix form as
\[
x_{it}'=\mu_{ix}+\gamma_{ix}'f_t+v_{itx}.
\]
Further let $\Gamma_i=(\lambda_i,\gamma_{ix})$,
$z_{it}=(y_{it},x_{it})'$, $\varepsilon_{it}=(e_{it},v_{itx}')'$, $\mu
_i=(\alpha_i,\mu_{ix}')'$. Then model (\ref{comshkmaineq1})
can be written as
\[
\lleft[\matrix{1 & -\beta'
\vspace*{3pt}\cr
0 & I_K}\rright]z_{it}=
\mu_i+\Gamma_i'f_t+
\varepsilon_{it}.
\]
Let $B$ denote the coefficient matrix of $z_{it}$ in the preceding
equation. Let $z_t=(z_{1t}',z_{2t}',\ldots,z_{Nt}')'$, $\Gamma
=(\Gamma
_1,\Gamma_2,\ldots,\Gamma_N)'$, $\varepsilon_t=(\varepsilon
_{1t}',\varepsilon_{2t}',\ldots,\varepsilon_{Nt}')'$ and $\mu=(\mu
_1',\mu_2',\ldots,\mu_N')'$. Stacking the equations over $i$, we have
\begin{equation}
\label{comshkmaineq2}(I_N\otimes B)z_t=\mu+\Gamma
f_t+\varepsilon_t.
\end{equation}
To analyze this model, we make the following assumptions.

\subsection{Assumptions}\label{sec2.1}

\renewcommand{\theass}{\Alph{ass}}
\begin{ass}\label{assA}
The factor process $f_t$ is a sequence of constants. Let
$M_{f
f}= T^{-1}\sum_{t=1}^T{\dot f_t\dot f_t'}$, where $\dot f_t=f_t-\frac
{1}{T}\sum_{t=1}^Tf_t$. We assume that $\widebar M_{f f} = \lim_{T \to\infty} M_{f f} $ is a strictly positive definite matrix.
\end{ass}

%
\begin{remark}
The nonrandomness assumption for $f_t$ is not
crucial. In fact, $f_t$ can be a sequence of random variables such that
$E(\|f_t\|^4)\le C<\infty$ uniformly in $t$, and $f_t$ is independent
of $\varepsilon_s$ for all $s$. The fixed $f_t$ assumption conforms
with the usual fixed effects assumption in panel data literature and,
in certain sense, is more general than random $f_t$.
\end{remark}

\begin{ass}\label{assB}
The idiosyncratic errors $\varepsilon
_{it}=(e_{it},v_{itx}')'$ are such that:
\begin{longlist}[(B.1)]
\item[(B.1)] The $e_{it}$ is independent and identically distributed over
$t$ and uncorrelated over $i$ with $E(e_{it})=0$ and $E(e_{it}^4)\le
\infty$ for all $i=1,\ldots, N$ and $t=1,\ldots, T$. Let $\Sigma_{iie}$
denote the variance of $e_{it}$.

\item[(B.2)] $v_{itx}$ is also independent and identically distributed
over $t$ and uncorrelated over $i$ with $E(v_{itx})=0$ and $E(\|
v_{itx}\|^4)\le\infty$ for all $i=1,\ldots, N$ and $t=1,\ldots, T$. We
use $\Sigma_{iix}$ to denote the variance matrix of $v_{itx}$.

\item[(B.3)] $e_{it}$ is independent of $v_{jsx}$ for all $(i,j,t,s)$.
Let $\Sigma_{ii}$ denote the variance matrix $\varepsilon_{it}$. So we
have $\Sigma_{ii}=\operatorname{diag}(\Sigma_{iie},\Sigma_{iix})$,
a block-diagonal matrix.
\end{longlist}
\end{ass}

%
\begin{remark}
Let $\Sigma_{\varepsilon\varepsilon}$ denote the
variance of $\varepsilon_t=(\varepsilon_{1t}',\ldots, \varepsilon
_{Nt}')'$. Due to the uncorrelatedness of $\varepsilon_{it}$ over $i$,
we have $\Sigma_{\varepsilon\varepsilon}=\operatorname{diag}(\Sigma
_{11},\Sigma
_{22},\ldots,\Sigma_{N N})$, a block-diagonal matrix. Assumption~\ref{assB} is
more general than the usual assumption in the factor analysis. In a
traditional factor model, the variances of the idiosyncratic error
terms are assumed to be a diagonal matrix. In the present setting, the
variance of $\varepsilon_t$ is a block-diagonal matrix. Even without
explanatory variables, this generalization is of interest. The factor
analysis literature has a long history to explore the block-diagonal
idiosyncratic variance, known as multiple battery factor analysis; see
\cite{tucker1958inter}. The maximum likelihood estimation theory for
high-dimensional factor models with block diagonal covariance matrix
has not been previously studied. The asymptotic theory developed in
this paper not only provides a way of analyzing the coefficient $\beta
$, but also a way of analyzing the factors and loadings in the multiple
battery factor models. This framework is of independent interest.
\end{remark}

\begin{ass}\label{assC}
There exists a $C>0$ sufficiently large such that:
\begin{longlist}[(C.3)]
\item[(C.1)]$\|\Gamma_j\| \le C$ for all $j=1,\ldots, N$;

\item[(C.2)]$C^{-1} \le\tau_{\min}(\Sigma_{j j}) \le\tau
_{\max}(\Sigma_{j j}) \le C$ for all $j=1,\ldots, N$, where $\tau
_{\min}(\Sigma_{j j})$ and $\tau_{\max}(\Sigma_{j j})$ denote the
smallest and largest eigenvalues of the matrix~$\Sigma_{j j}$, respectively;

\item[(C.3)] there exists an $r\times r$ positive matrix $Q$ such that
\[
Q=\lim_{N \to\infty} N^{-1}\Gamma'\Sigma_{\varepsilon
\varepsilon}^{-1}\Gamma,
\]
where $\Gamma$ is defined earlier.
\end{longlist}
\end{ass}

\begin{ass}\label{assD}
The variances $\Sigma_{ii}$ for all $i$ and
$M_{f f}$ are estimated in a compact set, that is, all the eigenvalues
of $\widehat\Sigma_{ii}$ and $\widehat M_{f f}$ are in an interval
$[C^{-1},C]$ for a sufficiently large constant $C$.
\end{ass}

\subsection{Identification restrictions}\label{sec2.2} It is a well-known result in
factor analysis that the factors and loadings can only be identified up
to a rotation; see, for example, \cite{anderson1956statistical,lawley1971factor}.
The models considered in this paper can be viewed as
extensions of the factor models. As such they inherit the same
identification problem. We show that identification conditions can be
imposed on the factors and loadings without loss of generality. To see
this, model (\ref{comshkmaineq2}) can be rewritten as
\begin{equation}
\label{add1} (I_N\otimes B)z_t=(\mu+\Gamma\bar f)+
\bigl[\Gamma M_{f f}^{1/2}R\bigr] \bigl[R'M_{f
f}^{-1/2}(f_t-
\bar f)\bigr]+\varepsilon_t,
\end{equation}
where $R$ is an orthogonal matrix, which we choose to be the matrix
consisting of the eigenvectors of $M_{f f}^{1/2}\Gamma'\Sigma
_{\varepsilon\varepsilon}^{-1}\Gamma M_{f f}^{1/2}$ associated with
the eigenvalues arranged in descending order. Treating $\mu+\Gamma
\bar
f$ as\vspace*{-2pt} the new $\mu^\star$, $\Gamma M_{f f}^{1/2}R$ as the new
$\Gamma
^\star$ and $R'M_{f f}^{-1/2}(f_t-\bar f)$ as the new $f_t^\star$,
we have
\[
(I_N\otimes B)z_t=\mu^\star+
\Gamma^\star f_t^\star+\varepsilon_t\vadjust{\goodbreak}
\]
with $\frac{1} T\sum_{t=1}^Tf_t^\star=0, \frac{1} T\sum_{t=1}^Tf_t^\star
f_t^{\star\prime}=I_r$ and $\frac{1} N\Gamma^{\star\prime}\Sigma
_{\varepsilon\varepsilon}^{-1}\Gamma^\star$ being a diagonal matrix.
Thus we impose the following restrictions for model (\ref
{comshkmaineq2}), which we refer to as \textup{IB} (\textit{identification}
restrictions for \textit{Basic} models).
\begin{longlist}[(IB1)]
\item[(IB1)] $M_{f f}=I_r$;

\item[(IB2)] $\frac{1}{N}\Gamma'\Sigma_{\varepsilon\varepsilon
}^{-1}\Gamma=D$, where $D$ is a diagonal matrix with its diagonal
elements distinct and arranged in descending order;

\item[(IB3)] $\bar f=\frac{1} T\sum_{t=1}^Tf_t=0$.
\end{longlist}

\subsection{Estimation}\label{sec2.3} The objective function considered in this
section is
\begin{equation}
\label{comshkmaineq3}\ln L(\theta)=-\frac{1}{2N}\ln |\Sigma_{z z} |-
\frac{1}{2N}\operatorname{tr} \bigl[(I_N\otimes
B)M_{z
z}\bigl(I_N\otimes B'\bigr)
\Sigma_{z z}^{-1} \bigr],
\end{equation}
where $\Sigma_{z z}=\Gamma M_{f f}\Gamma'+\Sigma_{\varepsilon
\varepsilon}$ and $M_{z z}=\frac{1} T\sum_{t=1}^T\dot z_t\dot z{}_t'$.
The latter is the data matrix. The parameters are $\theta=(\beta,
\Gamma, M_{f f}, \Sigma_{\varepsilon\varepsilon})$. The MLE is defined as
\[
\hat\theta=\mathop{\operatorname{argmax}}_{\theta\in\Theta}\ln L(\theta),
\]
where the parameter space $\Theta$ is defined to be a closed and
bounded subset containing the true parameter $\theta^*$ as an interior
point; $\Sigma_{\varepsilon\varepsilon}$ and $M_{f f}$ are positive
definite matrices, as in Assumption~\ref{assD}. The boundedness of $\Theta$
implies that the elements of $\beta$ and $\Gamma$ are bounded. This is
for theoretical purpose and is usually assumed for nonconvex
optimizations, as in \cite{jenrich1969asymptotic} and \cite
{newey1994large}. In actual computation with the EM algorithm, we do
not find the need to impose an upper or lower bound for the parameter
values. The likelihood function involves simple functions and are
continuous on $\Theta$ (in fact differentiable), so the MLE $\hat\theta
$ exists because a continuous function achieves its extreme value on a
closed and bounded subset.

Note that the determinant of $I_N\otimes B$ is 1, so the Jacobian term
does not depend on $B$. If $\varepsilon_t$ and $f_t$ are independent
and normally distributed, the likelihood function for the observed data
has the form of (\ref{comshkmaineq3}). Here recall that $f_t$ are fixed
constants, and $\varepsilon_t$ are not necessarily normal; (\ref
{comshkmaineq3}) is a pseudo-likelihood function.

For further analysis, we partition the matrix $\Sigma_{z z}$ and
$M_{z z}$ as\vspace*{6pt}
\[
\Sigma_{z z}=\pmatrix{
\Sigma_{z z}^{11} & \Sigma_{z z}^{12} & \cdots & \Sigma_{z z}^{1N}
\vspace*{3pt}\cr
\Sigma_{z z}^{21} & \Sigma_{z z}^{22} & \cdots& \Sigma_{z z}^{2N}
\vspace*{3pt}\cr
\vdots& \vdots& \ddots& \vdots
\vspace*{3pt}\cr
\Sigma_{z z}^{N1} & \Sigma_{z z}^{N2} & \cdots& \Sigma_{z z}^{N N}},
\qquad
M_{z z}=
\pmatrix{M_{z z}^{11} & M_{z z}^{12} & \cdots& M_{z z}^{1N}
\vspace*{3pt}\cr
M_{z z}^{21} & M_{z z}^{22} & \cdots&M_{z z}^{2N}
\vspace*{3pt}\cr
\vdots& \vdots& \ddots& \vdots
\vspace*{3pt}\cr
M_{z z}^{N1} & M_{z z}^{N2} & \cdots& M_{z z}^{N N}},
\]\vspace*{6pt}
where for any $(i,j)$, $\Sigma_{z z}^{i j}$ and $M_{z z}^{i j}$ are
both $(K+1)\times(K+1)$ matrices.\vadjust{\goodbreak}

Let $\hat\beta,\widehat\Gamma$ and $\widehat\Sigma_{\varepsilon
\varepsilon}$
denote the MLE. The first order condition for $\beta$ satisfies
\begin{equation}
\label{comshkmaineq4}
\frac{1}{NT}\sum_{i=1}^N
\sum_{t=1}^T\widehat\Sigma_{iie}^{-1}
\Biggl\{(\dot y_{it}- \dot x_{it}\hat\beta) -\hat
\lambda_i'\widehat G\sum_{j=1}^N
\widehat\Gamma_j\widehat\Sigma _{jj}^{-1}
\lleft[
\matrix{\dot y_{jt} -\dot x_{jt}\hat\beta
\vspace*{3pt}\cr
\dot x_{jt}'}\rright] \Biggr\}\dot x_{it}=0,\hspace*{-35pt}
\end{equation}
where $\widehat G=(\widehat M_{f f}^{-1}+\widehat\Gamma'\widehat\Sigma
_{\varepsilon
\varepsilon}^{-1}\widehat\Gamma)^{-1}$. The first order condition for
$\Gamma_j$ satisfies
\begin{equation}
\label{comshkmaineq5}\sum_{i=1}^N\widehat
\Gamma_i\widehat\Sigma _{ii}^{-1} \bigl(\widehat
BM_{z z}^{i j}\widehat B{}'-\widehat\Sigma_{z z}^{i
j}
\bigr)=0.
\end{equation}
Post-multiplying $\widehat\Sigma_{jj}^{-1}\widehat\Gamma_j'$ on both sides of
(\ref{comshkmaineq5}) and then taking summation over $j$, we have
\begin{equation}
\label{comshkmaineq6}\sum_{i=1}^N\sum
_{j=1}^N\widehat\Gamma _i\widehat\Sigma_{ii}^{-1} \bigl(\widehat BM_{z z}^{i j}
\widehat B{}'-\widehat\Sigma _{z
z}^{i j} \bigr)\widehat\Sigma_{j j}^{-1}\widehat\Gamma_j'=0.
\end{equation}
The first order condition for $\Sigma_{ii}$ satisfies
\begin{equation}
\label{comshkmaineq7} \widehat BM_{z z}^{ii}\widehat B{}'-\widehat\Sigma_{z z}^{ii}=\mathbb{W},
\end{equation}
where $\mathbb{W}$ is a $(K+1)\times(K+1)$ matrix such that its
upper-left $1\times1$ and lower-right $K\times K$ submatrices are both
zero, but the remaining elements are undetermined. The undetermined
elements correspond to the zero elements of $\Sigma_{ii}$. These first
order conditions are needed for the asymptotic representation of the MLE.

\subsection{Asymptotic properties of the MLE}\label{sec2.4} Theorem \ref{comshkthm2}
states the convergence rates of the MLE. The consistency is implied by
the theorem.

%
\begin{theorem}[(Convergence rate)]\label{comshkthm2}
Let $\hat\theta=(\hat\beta,\widehat\Gamma,\widehat\Sigma_{\varepsilon
\varepsilon
})$ be the solution by maximizing (\ref{comshkmaineq3}). Under
Assumptions~\ref{assA}--\ref{assD} and the identification conditions \textup{IB}, we have
\begin{eqnarray*}
\hat\beta-\beta&=&O_p\bigl(N^{-1/2}T^{-1/2}
\bigr)+O_p\bigl(T^{-1}\bigr),
\\
\frac{1}{N}\sum_{i=1}^N\bigl\|\widehat\Sigma_{ii}^{-1}\bigr\|\cdot\|\widehat\Gamma _i-
\Gamma_i\|^2&=&O_p\bigl(T^{-1}
\bigr),\qquad
\frac{1}{N}\sum_{i=1}^N\|\widehat\Sigma_{ii}-\Sigma_{ii}\|^2 =O_p
\bigl(T^{-1}\bigr).
\end{eqnarray*}
\end{theorem}

%
\begin{remark}
Bai \cite{bai2009panel} considers an iterated
principal components estimator for model (\ref{comshkmaineq1}). His
derivation shows that, in the presence of heteroscedasticities over the
cross section, the PC estimator for $\beta$ has a bias of order
$O_p(N^{-1})$. As a~comparison, Theorem \ref{comshkthm2} shows that the
MLE is robust to the heteroscedasticities over the cross section. So if
$N$ is fixed, the estimator in \cite{bai2009panel} is inconsistent
unless there is no heteroskedasticity, but the estimator here is still
\mbox{consistent}.
\end{remark}

Let $\mathcal{M}(\mathbb X)$ denote the project matrix onto the space
orthogonal to $\mathbb X$, that is, $\mathcal{M}(\mathbb X)=I-\mathbb
X(\mathbb X'\mathbb X)^{-1}\mathbb X'$. We have

%
\begin{theorem}[(Asymptotic representation)]\label{comshkthm3}
Under the assumptions of Theorem \ref{comshkthm2}, we have
\begin{eqnarray*}
\hat\beta-\beta&=&\Omega^{-1}\frac{1}{NT}\sum
_{i=1}^N\sum_{t=1}^T
\Sigma _{iie}^{-1}e_{it}v_{itx}
\\
&&{} +O_p\bigl(T^{-3/2}\bigr)+O_p\bigl(N^{-1}T^{-1/2}
\bigr)+O_p\bigl(N^{-1/2}T^{-1}\bigr),
\end{eqnarray*}
where $\Omega$ is a $K\times K$ matrix whose $(p,q)$ element $\Omega
_{pq}=\frac{1}{N}\sum_{i=1}^N\Sigma_{iie}^{-1}\Sigma_{iix}^{(p,q)}$
with $\Sigma_{iix}^{(p,q)}$ being the $(p,q)$ element of matrix
$\Sigma_{iix}$.
\end{theorem}

%
\begin{remark}\label{rmk25}
In Appendix~A.3 of the supplement~\cite{bai}, we show that the asymptotic expression of $\hat\beta-\beta$ can
be alternatively expressed as
\begin{eqnarray}\label{alternative1}
\hat\beta-\beta&=&\pmatrix{
\operatorname{tr}\bigl[\ddot MX_1\mathcal {M}(\widebar{\mathbb F})
X_1' \bigr]& \cdots& \operatorname{tr}\bigl[\ddot
MX_1\mathcal{M}(\widebar{\mathbb F}) X_K'
\bigr]
\vspace*{3pt}\cr
\vdots& \vdots& \vdots
\vspace*{3pt}\cr
\operatorname{tr}\bigl[\ddot MX _K\mathcal {M}(\widebar{\mathbb F})X_1' \bigr] & \cdots&\operatorname{tr}\bigl[\ddot MX
_K\mathcal {M}(\widebar{\mathbb F}) X _K'
\bigr]}^{-1}
\nonumber\\[20pt]\\[-30pt]
&&{} \times\pmatrix{
\operatorname{tr}\bigl[\ddot
MX_1\mathcal{M}(\widebar{\mathbb F}) e' \bigr]
\vspace*{3pt}\cr
\vdots
\vspace*{3pt}\cr
\operatorname{tr}\bigl[\ddot MX _K\mathcal{M}(\widebar{\mathbb F}) e'\bigr]}\nonumber
\\
&&{} +O_p \bigl(T^{-3/2}\bigr)+O_p\bigl(N^{-1}T^{-1/2}\bigr)+O_p\bigl(N^{-1/2}T^{-1}\bigr),\nonumber
\end{eqnarray}
where $X_k =  (x_{itk} )$ is $N\times T$ (the data matrix for
the $k$th regressor, $k=1,2,\ldots, K$); $e= (e_{it} )$ is
$N\times T$; $\ddot M=\Sigma_{ee}^{-1/2}\mathcal{M}(\Sigma
_{ee}^{-1/2}\Lambda)\Sigma_{ee}^{-1/2}$ with $\Sigma
_{ee}=\operatorname{diag}\{\Sigma
_{11e}, \Sigma_{22e},\ldots,\Sigma_{N Ne}\}$ and $\Lambda
=(\lambda
_1, \lambda_2, \ldots, \lambda_N)'$; $\mathbb F=(f_1,f_2,\ldots, f_T)'$;
$\widebar{\mathbb F}=(1_T, \mathbb F)$ where $1_T$ is a $T\times1$
vector with all 1's.
\end{remark}

%
\begin{remark}
Theorem \ref{comshkthm3} shows that the asymptotic
expression of $\hat\beta-\beta$ only involves variations in $e_{it}$
and $v_{itx}$. Intuitively, this is due to the fact that the error
terms of the $y$ equation share the same factors with the explanatory
variables. The variations from the common factor part of $x_{itk}$
(i.e., $\gamma_{ik}'f_t$) do not provide information for $\beta$ since
this part of information is offset by the common factor part of the
error terms (i.e., $\lambda_i'f_t$) in the $y$ equation.
\end{remark}

%
\begin{corollary}[(Limiting distribution)]\label{comshkthm4}
Under the assumptions of Theorem~\ref{comshkthm3}, if $\sqrt N/T\to0$,
we have
\[
\sqrt{NT} (\hat\beta-\beta)\stackrel{d} {\rightarrow }N\bigl(0,{\widebar\Omega}^{ -1}\bigr),
\]
where $\widebar{\Omega}=
\lim_{N,T\to\infty}\Omega$, and $\widebar\Omega$ is also the
limit of
\[
\widebar{\Omega}=\mathop{\operatorname{plim}}_{N,T\to\infty
}\frac{1} {NT}
\pmatrix{
\operatorname{tr}\bigl[\ddot
MX_1\mathcal {M}(\widebar{\mathbb F}) X_1'
\bigr]& \cdots& \operatorname{tr}\bigl[\ddot MX_1\mathcal{M}(
\widebar{\mathbb F}) X_K' \bigr]
\vspace*{3pt}\cr
\vdots& \vdots& \vdots
\vspace*{3pt}\cr
\operatorname{tr}\bigl[\ddot MX _K\mathcal {M}(\widebar{\mathbb F})X_1' \bigr] & \cdots&\operatorname{tr}\bigl[\ddot MX
_K\mathcal {M}(\widebar{\mathbb F}) X _K'
\bigr]}.
\]
\end{corollary}

%
\begin{remark}\label{remark2}
 Matrix $\widebar\Omega$ can be
consistently estimated by
\[
\frac{1} {NT}\pmatrix{
\operatorname{tr}\bigl[\widehat{\ddot M}X_1\mathcal {M}(\widehat{\widebar{\mathbb F}}) X_1' \bigr]& \cdots&
\operatorname{tr}\bigl[\widehat{\ddot M}X_1\mathcal{M}(\widehat{\widebar{\mathbb F}}) X_K' \bigr]
\vspace*{3pt}\cr
\vdots& \vdots& \vdots
\vspace*{3pt}\cr
\operatorname{tr}\bigl[\widehat{\ddot M}X_K\mathcal{M}(\widehat{\widebar{\mathbb F}})X_1' \bigr] & \cdots&
\operatorname{tr}\bigl[\widehat{\ddot M}X_K\mathcal {M}(\widehat{\widebar{\mathbb F}}) X _K' \bigr]},
\]
where $X_k$ is the $N\times T$ data matrix for the $k$th regressor,
\begin{equation}
\label{1110a}\widehat{\ddot M}=\widehat\Sigma _{ee}^{-1}-\widehat\Sigma_{ee}^{-1}\widehat\Lambda \bigl(\widehat\Lambda'
\widehat\Sigma_{ee}^{-1}\widehat\Lambda\bigr)^{-1}\widehat\Lambda '\widehat\Sigma _{ee}^{-1};
\end{equation}
$\widehat{\widebar{\mathbb F}}=(1_T, \widehat{\mathbb F})$ with $\widehat
{\mathbb F}=(\hat f_1, \hat f_2, \ldots, \hat f_T)'$ and
\begin{equation}
\label{1110b}\hat f_t=\Biggl(\sum_{i=1}^N
\widehat\Gamma_i\widehat\Sigma_{ii}^{-1}\widehat\Gamma_i'\Biggr)^{-1} \Biggl(\sum
_{i=1}^N\widehat\Gamma_i\widehat\Sigma_{ii}^{-1}\widehat B\dot z_{it}\Biggr).
\end{equation}
Here $\widehat\Gamma, \widehat\Lambda, \widehat\Sigma_{ii}, \widehat\Sigma
_{ee}$ and
$\widehat B$ are the maximum likelihood estimators.
\end{remark}

\section{Common shock models with zero restrictions}\label{sec3}
The basic model in Section~\ref{sec2} assumes that the explanatory variables
$x_{it}$ share the same factors with $y_{it}$. This section relaxes
this assumption. We assume that the regressors are impacted by
additional factors that do not affect the $y$ equation. An alternative
view is that some factor loadings in the $y$ equation are restricted to
be zero. Consider the following model:
\begin{eqnarray}
\label{zerores1} y_{it}&=&\alpha_i+x_{it1}
\beta_1+x_{it2}\beta_2+\cdots+x_{itK}
\beta _K+\psi_i'g_t+e_{it},
\nonumber
\\[-8pt]
\\[-8pt]
x_{itk}&=&\mu_{ik}+\gamma_{ik}^{g\prime}g_t+
\gamma_{ik}^{h\prime
}h_t+v_{itk}
\nonumber
\end{eqnarray}
for $k=1,2,\ldots,K$, where $g_t$ is an $r_1\times1$ vector
representing the shocks affecting both $y_{it}$ and $x_{it}$, and $h_t$
is an $r_2\times1$ vector representing the shocks affecting $x_{it}$
only. Let $\lambda_i=(\psi_i',0_{r_2\times1}')'$, $\gamma
_{ik}=(\gamma
_{ik}^{g\prime},\gamma_{ik}^{h\prime})'$ and $f_t=(g_t',h_t')'$, the
above model can be written as
\begin{eqnarray*}
y_{it}&=&\alpha_i+x_{it1}\beta_1+x_{it2}
\beta_2+\cdots +x_{itK}\beta_K+
\lambda_i'f_t+e_{it},
\\
x_{itk}&=&\mu_{ik}+\gamma_{ik}'f_t+v_{itk},
\end{eqnarray*}
which is the same as model (\ref{comshkmaineq1}) except that $r_2$
elements of $\lambda_i$ are restricted to be zeros. For further
analysis, we introduce some notation. We define
\begin{eqnarray*}
\Gamma_i^g&=&\bigl(\psi_i,
\gamma_{i1}^g, \ldots, \gamma_{iK}^g
\bigr),\qquad \Gamma_i^h=\bigl(0_{r_2\times1},
\gamma_{i1}^h, \ldots,\gamma _{iK}^h
\bigr),
\\
\Gamma^g&=&\bigl(\Gamma_1^g,
\Gamma_2^g, \ldots, \Gamma_N^g
\bigr)',\qquad \Gamma^h=\bigl(\Gamma_1^h,
\Gamma_2^h,\ldots, \Gamma_N^h
\bigr)'.
\end{eqnarray*}
We also define $\mathbb G$ and $\mathbb H$ similarly as $\mathbb F$,
that is, $\mathbb G=(g_1, g_2,\ldots, g_T)'$, $\mathbb H=(h_1, h_2,\ldots, h_T)'$. This implies that $\mathbb F=(\mathbb G, \mathbb H)$.
The presence of zero restrictions in (\ref{zerores1}) requires
different identification conditions.

\subsection{Identification conditions}\label{sec3.1} Zero loading restrictions
alleviate rotational indeterminacy. Instead of $r^2=(r_1+r_2)^2$
restrictions, we only need to impose $r_1^2+r_1r_2+r_2^2$ restrictions.
These restrictions are referred to as IZ restrictions (\textit{Identification} conditions with \textit{Zero} restrictions). They are:
\begin{longlist}[(IZ2)]
\item[(IZ1)] $M_{f f}=I_r$;
\item[(IZ2)] $\frac{1}{N}\Gamma^{g\prime}\Sigma_{\varepsilon
\varepsilon
}^{-1}\Gamma^g=D_1$ and $\frac{1}{N}\Gamma^{h\prime}\Sigma
_{\varepsilon
\varepsilon}^{-1}\Gamma^h=D_2$, where $D_1$ and $D_2$ are both diagonal
matrices with distinct diagonal elements in descending order;
\item[(IZ3)] $1_T'\mathbb G=0$ and $1_T'\mathbb H=0$.\vadjust{\goodbreak}
\end{longlist}

In addition, we need an additional assumption for our analysis.

\begin{ass}\label{assE} $\Psi=(\psi_1',\psi_2',\ldots,\psi_N')'$ is of
full column rank.
\end{ass}

Identification conditions IZ are less stringent than \textup{IB} of the previous
section. Assumption~\ref{assE} says that the factors $g_t$ are pervasive for the
$y$ equation. In Appendix~B of the supplement~\cite{bai}, we explain why
$r_1^2+r_1r_2+r_2^2$ restrictions are sufficient.

\subsection{Estimation}\label{sec3.2}
The likelihood function is now maximized under three sets of
restrictions, that is, $\frac{1}{N}\Gamma^{g\prime}\Sigma
_{\varepsilon
\varepsilon}^{-1}\Gamma^g=D_1$, $\frac{1}{N}\Gamma^{h\prime}\Sigma
_{\varepsilon\varepsilon}^{-1}\Gamma^h=D_2$ and $\Phi=0$ where
$\Phi$
denotes the zero factor loading matrix in the $y$ equation.
The likelihood function with the Lagrange multipliers is
\begin{eqnarray*}
\ln L&=&-\frac{1}{2N}\ln|\Sigma_{z z}|-\frac{1}{2N}
\operatorname {tr} \bigl[(I_N\otimes B)M_{z z}
\bigl(I_N\otimes B'\bigr)\Sigma_{z z}^{-1}
\bigr]
\\
&&{} +\operatorname{tr} \biggl[\Upsilon_1 \biggl(\frac{1}{N}
\Gamma^{g\prime
}\Sigma_{\varepsilon
\varepsilon}^{-1}\Gamma^g
-D_1 \biggr) \biggr]+\operatorname{tr} \biggl[\Upsilon_2
\biggl(\frac
{1}{N}\Gamma^{h\prime}\Sigma _{\varepsilon\varepsilon}^{-1}
\Gamma^h-D_2 \biggr) \biggr]
\\
&&{} +\operatorname {tr}\bigl[
\Upsilon_3'\Phi\bigr],
\end{eqnarray*}
where $\Sigma_{z z}=\Gamma\Gamma'+\Sigma_{\varepsilon\varepsilon}$;
$\Upsilon_1$ is $r_1\times r_1$ and $\Upsilon_2$ is $r_2\times r_2$,
both are symmetric Lagrange multipliers matrices with zero diagonal
elements; $\Upsilon_3$ is a Lagrange multiplier matrix of dimension
$r_2\times N$.

Let $\mathbb U=\widehat\Sigma_{z z}^{-1}[(I_N\otimes\widehat B)M_{z
z}(I_N\otimes\widehat B{}')-\widehat\Sigma_{z z}]\widehat\Sigma_{z z}^{-1}$.
Notice $\mathbb U$ is a symmetric matrix. The first order condition on
$\widehat\Gamma^g$ gives
\[
\frac{1}{N}\widehat\Gamma^{g\prime}\mathbb U+\Upsilon_1
\frac
{1}{N}\widehat\Gamma ^{g\prime}\widehat\Sigma_{\varepsilon\varepsilon}^{-1}=0.
\]
Post-multiplying $\widehat\Gamma^g$ yields
\[
\frac{1}{N}\widehat\Gamma^{g\prime}\mathbb U\widehat\Gamma^g+
\Upsilon _1\frac
{1}{N}\widehat\Gamma^{g\prime}\widehat\Sigma_{\varepsilon\varepsilon
}^{-1}\widehat\Gamma^g=0.
\]
Since $\frac{1}{N}\widehat\Gamma^{g\prime}\mathbb U\widehat\Gamma^g$ is a
symmetric matrix, the above equation implies that $\Upsilon_1\frac
{1}{N}\widehat\Gamma^{g\prime}\widehat\Sigma_{\varepsilon\varepsilon
}^{-1}\widehat\Gamma^g$ is also symmetric. But $\frac{1}{N}\widehat\Gamma^{g\prime
}\widehat\Sigma_{\varepsilon\varepsilon}^{-1}\widehat\Gamma^g$ is a diagonal matrix.
So the $(i,j)$th element of $\Upsilon_1\frac{1}{N}\widehat\Gamma
^{g\prime
}\widehat\Sigma_{\varepsilon\varepsilon}^{-1}\widehat\Gamma^g$ is
$\Upsilon
_{1,i j}d_{1j}$, where $\Upsilon_{1,i j}$ is the $(i,j)$th element of
$\Upsilon_1$ and $d_{1j}$ is the $j$th diagonal element of $\widehat D_1$.
Given $\Upsilon_1\frac{1}{N}\widehat\Gamma^{g\prime}\widehat\Sigma
_{\varepsilon
\varepsilon}^{-1}\widehat\Gamma^g$ is symmetric, we have $\Upsilon
_{1,i
j}d_{1j}=\Upsilon_{1,ji}d_{1i}$ for all $i\neq j$. However, $\Upsilon
_1$ is also symmetric, so $\Upsilon_{1,i j}=\Upsilon_{1,ji}$. This
gives $\Upsilon_{1,i j}(d_{1j}-d_{1i})=0$. Since $d_{1j}\neq d_{1i}$
by IZ2, we have $\Upsilon_{1,i j}=0$ for all $i\neq j$. This
implies $\Upsilon_1=0$ since the diagonal elements of $\Upsilon_1$ are
all zeros.

Let $\Gamma_x^h=(\gamma_{1x}^h,\gamma_{2x}^h,\ldots,\gamma_{Nx}^h)'$
with $\gamma_{ix}^h=(\gamma_{i1}^h,\gamma_{i2}^h,\ldots,\gamma
_{iK}^h)$, and $\Sigma_{xx}=\break \operatorname{diag}\{\Sigma_{11x},\Sigma
_{22x},\ldots,\Sigma_{N Nx}\}$, a block diagonal matrix of $NK\times NK$ dimension.
We partition the matrix $\mathbb U$ and define the matrix $\widebar{\mathbb{U}}$ as
\[
\mathbb U=\pmatrix{
\mathbb U_{11} &\mathbb U_{12} & \cdots&\mathbb U_{1N}
\vspace*{3pt}\cr
\mathbb U_{21} &\mathbb U_{22} & \cdots&\mathbb U_{2N}
\vspace*{3pt}\cr
\vdots& \vdots& \ddots& \vdots
\vspace*{3pt}\cr
\mathbb U_{N1} &\mathbb U_{N2} & \cdots&\mathbb
U_{N N}},
\qquad
\widebar{\mathbb{U}}=
\pmatrix{\widebar{\mathbb{U}}_{11} &\widebar{\mathbb U}_{12} & \cdots& \widebar{\mathbb U}_{1N}
\vspace*{3pt}\cr
\widebar{\mathbb U}_{21} & \widebar{\mathbb U}_{22} & \cdots& \widebar{\mathbb U}_{2N}
\vspace*{3pt}\cr
\vdots& \vdots& \ddots& \vdots
\vspace*{3pt}\cr
\widebar{\mathbb U}_{N1} & \widebar{\mathbb U}_{N2} &
\cdots& \widebar{\mathbb U}_{N N}},
\]
where $\mathbb U_{i j}$ is a $(K+1)\times(K+1)$ matrix, and
$\widebar{\mathbb U}_{i j}$ is the lower-right $K\times K$ block of $\mathbb
U_{i j}$. Notice $\widebar{\mathbb U}$ is also a symmetric matrix.
Then the first order condition on $\Gamma_x^h$ gives
\[
\frac{1}{N}\widehat\Gamma_x^{h\prime}\widebar{\mathbb U}+
\Upsilon _2\frac
{1}{N}\widehat\Gamma_x^{h\prime}
\widehat\Sigma_{xx}^{-1}=0.
\]
Post-multiplying $\widehat\Gamma_x^h$ yields
\[
\frac{1}{N}\widehat\Gamma_x^{h\prime}\widebar{\mathbb U}
\widehat\Gamma _x^h+\Upsilon_2\frac{1}{N}
\widehat\Gamma_x^{h\prime}\widehat\Sigma _{xx}^{-1}
\widehat\Gamma_x^h=0.
\]
Notice $\frac{1}{N}\widehat\Gamma_x^{h\prime}\widehat\Sigma_{xx}^{-1}\widehat\Gamma
_x^h=\frac{1}{N}\widehat\Gamma^{h\prime}\widehat\Sigma_{\varepsilon
\varepsilon
}^{-1}\widehat\Gamma^h=\widehat D_2$. By the similar arguments in deriving
$\Upsilon_1=0$, we have $\Upsilon_2=0$. The interpretation for the zero
Lagrange multipliers is that these constraints do not affect the
optimal value of the likelihood function nor the efficiency of $\hat
\beta$. In contrast, we cannot show $\Upsilon_3$ to be zero. Thus the
restriction $\Phi=0$ affects the optimal value of the likelihood
function and the efficiency of $\hat\beta$. In Section~\ref{sec2}, we did not
use the Lagrange multiplier approach to analyze the identification
restrictions. Had this been done, we would have obtained zero valued
Lagrange multipliers. This is another view of why these restrictions do
not affect the limiting distribution of $\hat\beta$. But these
restrictions are needed to remove the rotational indeterminacy. 

Now the likelihood function is simplified as
\begin{equation}
\label{extenmainadd}\qquad \ln L=-\frac{1}{2N}\ln|\Sigma _{z
z}|-
\frac{1}{2N}\operatorname{tr} \bigl[(I_N\otimes
B)M_{z
z}\bigl(I_N\otimes B'\bigr)
\Sigma_{z
z}^{-1} \bigr]+\operatorname{tr}\bigl[
\Upsilon_3'\Phi\bigr].
\end{equation}
The first order condition on $\Gamma$ is
\begin{equation}
\label{extenmaineq1}\widehat\Gamma{}'\widehat\Sigma_{z
z}^{-1}
\bigl[(I_N\otimes\widehat B)M_{z z}\bigl(I_N\otimes
\widehat B{}'\bigr)-\widehat\Sigma _{z
z}\bigr]\widehat\Sigma_{z z}^{-1}=W',
\end{equation}
where $W$ is a matrix having the same dimension as $\Gamma$, whose
element is zero if the counterpart of $\Gamma$ is not specified to be
zero, otherwise undetermined (containing the Lagrange multipliers).
Post-multiplying $\widehat\Gamma$ gives
\[
\widehat\Gamma'\widehat\Sigma_{z z}^{-1}
\bigl[(I_N\otimes\widehat B)M_{z
z}\bigl(I_N\otimes
\widehat B{}'\bigr)-\widehat\Sigma_{z z}\bigr]\widehat\Sigma_{z z}^{-1}\widehat\Gamma=W'\widehat\Gamma.
\]
By the special structure of $W$ and $\widehat\Gamma$, it is easy to verify
that $W'\widehat\Gamma$ has the form
\[
\lleft[\matrix{0_{r_1\times r_1} & 0_{r_1\times r_2}
\vspace*{3pt}\cr
\times&
0_{r_2\times r_2}} \rright].
\]
However, the left-hand side of the preceding equation is a symmetric
matrix, and so is the right-hand side. It follows that the subblock
``$\times$'' is zero, that is, \mbox{$W'\widehat\Gamma=0$}.\vadjust{\goodbreak} Thus,
$\widehat\Gamma'\widehat\Sigma_{z z}^{-1}[(I_N\otimes\widehat B)M_{z
z}(I_N\otimes\widehat B{}')-\widehat\Sigma_{z z}]\widehat\Sigma_{z z}^{-1}\widehat\Gamma=0$. (This equation would be the first order condition for
$M_{f
f}$ if it were unknown.) This equality can be simplified as
\begin{equation}
\label{extenmaineq2}\widehat\Gamma'\widehat\Sigma _{\varepsilon
\varepsilon}^{-1}
\bigl[(I_N\otimes\widehat B)M_{z z}\bigl(I_N\otimes
\widehat B{}'\bigr)-\widehat\Sigma_{z z}\bigr]\widehat\Sigma_{\varepsilon\varepsilon}^{-1}\widehat\Gamma=0
\end{equation}
because $\widehat\Gamma'\widehat\Sigma_{z z}^{-1}=\widehat G\widehat\Gamma'\widehat\Sigma
_{\varepsilon\varepsilon}^{-1}$ with $\widehat G=(I+\widehat\Gamma'\widehat\Sigma
_{\varepsilon\varepsilon}^{-1}\widehat\Gamma)^{-1}$. Next, we partition the
matrix $\widehat G=(I+\widehat\Gamma'\widehat\Sigma_{\varepsilon\varepsilon
}^{-1}\widehat\Gamma)^{-1}$ and $\widehat H=(\widehat\Gamma'\widehat\Sigma
_{\varepsilon
\varepsilon}^{-1}\widehat\Gamma)^{-1}$ as follows:
\[
\widehat G= \lleft[\matrix{\widehat G_1
\vspace*{3pt}\cr
\widehat G_2} \rright] =
\lleft[\matrix{\widehat G_{11} & \widehat G_{12}
\vspace*{3pt}\cr
\widehat G_{21} & \widehat G_{22}} \rright], \qquad \widehat H= \lleft[
\matrix{\widehat H_1
\vspace*{3pt}\cr
\widehat H_2} \rright] = \lleft[\matrix{
\widehat H_{11} & \widehat H_{12}
\vspace*{3pt}\cr
\widehat H_{21} & \widehat H_{22}} \rright],
\]
where $\widehat G_{11}, \widehat H_{11}$ are $r_1\times r_1$, while $\widehat G_{22}, \widehat H_{22}$ are $r_2\times r_2$.

Notice $\widehat\Sigma_{z z}^{-1}=\widehat\Sigma_{\varepsilon\varepsilon
}^{-1}-\widehat\Sigma_{\varepsilon\varepsilon}^{-1}\widehat\Gamma\widehat G\widehat\Gamma'\widehat\Sigma_{\varepsilon\varepsilon}^{-1}$ and $\widehat\Gamma
'\widehat\Sigma_{z z}^{-1}=\widehat G\widehat\Gamma'\widehat\Sigma_{\varepsilon
\varepsilon
}^{-1}$. Substitute these results into (\ref{extenmaineq1}), and use
(\ref{extenmaineq2}). The first order condition for $\psi_i$ can be
simplified as
\begin{equation}
\label{extenmaineq3}\widehat G_{1}\sum_{i=1}^N
\widehat\Gamma _i\widehat\Sigma_{ii}^{-1}\bigl(\widehat BM_{z z}^{i j}\widehat B{}'-\widehat\Sigma_{z
z}^{i j}
\bigr)\widehat\Sigma_{j j}^{-1}I_{K+1}^1=0,
\end{equation}
where $I _{K+1}^{ 1}$ is the first column of the identity matrix of
dimension $K+1$.

Similarly, the first order condition for $\gamma_{jx}=(\gamma
_{j1},\gamma_{j2},\ldots,\gamma_{jK})$ is
\begin{equation}
\label{extenmaineq4}\sum_{i=1}^N\widehat\Gamma_i\widehat\Sigma _{ii}^{-1}\bigl(\widehat BM_{z z}^{i j}\widehat B{}'-\widehat\Sigma_{z z}^{i j}
\bigr)\widehat\Sigma_{j j}^{-1}I _{K+1}^{ -}=0,
\end{equation}
where $I _{K+1}^{ -}$ is a $(K+1)\times K$ matrix, obtained by
deleting the first column of the identity matrix of dimension $K+1$.

The first order condition for $\Sigma_{j j}$ is
\begin{eqnarray}
\label{extenmaineq5} && \widehat B M_{z z}^{j j}\widehat B{}'-\widehat\Sigma_{z z}^{j
j}-\widehat\Gamma_j'\widehat G\sum_{i=1}^N \widehat\Gamma_i
\widehat\Sigma _{ii}^{-1} \bigl(\widehat B M_{z z}^{i j}
\widehat B{}'-\widehat\Sigma_{z z}^{i j} \bigr)
\nonumber
\\[-8pt]
\\[-8pt]
&&\qquad{} -\sum_{i=1}^N \bigl(\widehat BM_{z z}^{ji}\widehat B{}'-\widehat\Sigma_{z
z}^{ji}
\bigr)\widehat\Sigma_{ii}^{-1}\widehat\Gamma_i'
\widehat G\widehat\Gamma_j=\mathbb {W},
\nonumber
\end{eqnarray}
where $\mathbb{W}$ is defined following (\ref{comshkmaineq7}).

The first order condition for $\beta$ is
\begin{equation}
\label{extenmainbeta}
\frac{1}{NT}\sum_{i=1}^N
\sum_{t=1}^T\widehat\Sigma_{iie}^{-1}
\Biggl\{(\dot y_{it}- \dot x_{it}\hat\beta) -\hat
\lambda_i'\widehat G\sum_{j=1}^N
\widehat\Gamma_j\widehat\Sigma _{jj}^{-1}
\lleft[
\matrix{\dot y_{jt} -\dot x_{jt}\hat
\beta
\vspace*{3pt}\cr
\dot x_{jt}'}\rright]
\Biggr\}\dot x_{it}=0,\hspace*{-35pt}
\end{equation}
which is the same as in Section~\ref{sec2}.

We need an additional identity to study the properties of the MLE.
Recall that, by the special structures of $W$ and $\widehat\Gamma$, the
three submatrices of $W'\widehat\Gamma$ can be directly derived to be
zeros. The remaining submatrix is also zero, as shown earlier. However,
this submatrix being zero yields the following equation (the detailed
derivation is delivered in Appendix~F):
\begin{equation}
\label{extenmaineq6}\frac{1}{N}\widehat G_{2}\sum
_{i=1}^N\sum_{j=1}^N
\widehat\Gamma_i\widehat\Sigma_{ii}^{-1}\bigl(\widehat BM_{z
z}^{i
j}\widehat B{}'-\widehat\Sigma_{z z}^{i j}
\bigr)\widehat\Sigma_{j j}^{-1}I _{K+1}^{
1}
\hat\psi_j'=0.
\end{equation}
These identities are used to derive the asymptotic representations.

\subsection{Asymptotic properties of the MLE}\label{sec3.3} The results on
consistency and the rate of convergence are similar to those in the
previous section, which are presented in Appendixes B.1 and B.2. For
simplicity, we only state the asymptotic representation for the MLE here.

%
\begin{proposition}[(Asymptotic representation)]\label{extenthm1}
$\!\!\!\!$Under Assumptions \mbox{\ref{assA}--\ref{assE}} and the identification restriction IZ, we have
\begin{eqnarray*}
\mathcal{P}^0(\hat\beta-\beta)&=&\frac{1}{NT}\sum
_{i=1}^N\sum_{t=1}^T
\Sigma_{iie}^{-1}e_{it}v_{itx}+
\frac{1}{NT} \sum_{i=1}^N\sum
_{t=1}^T\Sigma_{iie}^{-1}
\gamma_{ix}^{h\prime} h_te_{it}
\\
&&{}-\frac{1}{NT}\sum_{i=1}^N\sum
_{t=1}^T\Sigma_{iie}^{-1}
\psi _i'\Pi_{\psi
\psi}^{-1} \Biggl(
\frac{1}{N}\sum_{j=1}^N
\psi_j\Sigma_{j
je}^{-1}\gamma
_{jx}^{h\prime} \Biggr) h_te_{it}
\\
&&{}+O_p\bigl(T^{-3/2}\bigr)+O_p
\bigl(N^{-1}T^{-1/2}\bigr)+O_p\bigl(N^{-1/2}T^{-1}
\bigr),
\end{eqnarray*}
where $\mathcal{P}^0$ is a $K\times K$ symmetric matrix with its
$(p,q)$ element equal to $\frac{1} N\operatorname{tr}(\Gamma
_p^{h\prime}\ddot M\Gamma
_q^h)+\frac{1} N\sum_{i=1}^N\Sigma_{iie}^{-1}\Sigma_{iix}^{(p,q)}$;
$\Gamma_p^h=[\gamma_{1p}^h,\gamma_{2p}^h,\ldots, \gamma_{Np}^h]'$;
$\gamma_{jx}^h=[\gamma_{j1}^h,\ldots,\break  \gamma_{jK}^h]$; $\Pi_{\psi
\psi
}=\frac{1}{N}\sum_{i=1}^N\psi_i
\Sigma_{iie}^{-1}\psi_i'$ and $\ddot M=\Sigma_{ee}^{-1/2}\mathcal
M(\Sigma_{ee}^{-1/2}\Psi)\Sigma_{ee}^{-1/2}$.
\end{proposition}

Proposition \ref{extenthm1} is derived under the identification
conditions IZ. In Appendix~B.3 of the supplement~\cite{bai},
we show that for any set of factors and factor loadings $(\psi_i,
\gamma
_{ik}, g_t, h_t)$, it can always be transformed into a new set $(\psi
_i^\star, \gamma_{ik}^\star, g_t^\star, h_t^\star)$, which satisfies
IZ, and at the same time, leaving $\Phi=0$ intact. Given the asymptotic
representation in Proposition \ref{extenthm1}, together with the
relationship between the two sets, we have the following theorem, which
does not depend on IZ.

%
\begin{theorem}\label{extenthmnew}
Under Assumptions \ref{assA}--\ref{assE}, we have
\begin{eqnarray*}
\mathcal{P}(\hat\beta-\beta)&=& \frac{1}{NT}\sum
_{i=1}^N\sum_{t=1}^T
\Sigma_{iie}^{-1}e_{it}v_{itx} +
\frac{1}{NT}\sum_{i=1}^N\sum
_{t=1}^T\Sigma_{iie}^{-1}\gamma
_{ix}^{h\prime}h_t^\star e_{it}
\\
&&{}-\frac{1}{NT}\sum_{i=1}^N\sum
_{t=1}^T\Sigma_{iie}^{-1}
\psi _i'\Pi_{\psi
\psi}^{-1} \Biggl(
\frac{1}{N}\sum_{j=1}^N
\psi_j\Sigma_{j je}^{-1}\gamma
_{jx}^{h\prime} \Biggr)h_t^\star
e_{it}
\\
&&{}+O_p\bigl(T^{-3/2}\bigr)+O_p
\bigl(N^{-1}T^{-1/2}\bigr)+O_p\bigl(N^{-1/2}T^{-1}
\bigr),
\end{eqnarray*}
where
\[
h_t^\star=\dot h_t-\dot{\mathbb
H}'\dot{\mathbb G}\bigl(\dot{\mathbb G}'\dot {\mathbb
G}\bigr)^{-1}\dot g_t;
\]
$\mathcal{P}$ is a $K\times K$ symmetric matrix with its $(p,q)$
element equal to
\[
\frac{1} {NT}\operatorname{tr} \bigl[\ddot M\Gamma_q^h\mathbb
H'\mathcal M(\widebar{\mathbb G})\mathbb H\Gamma_p^{h\prime}
\bigr]+\frac{1} N\sum_{i=1}^N\Sigma
_{iie}^{-1}\Sigma_{iix}^{(p,q)},
\]
where $\widebar{\mathbb G}=(1_T, \mathbb G)$; $\Pi_{\psi\psi
}=\frac{1} N
\sum_{i=1}^N\psi_i\Sigma_{iie}^{-1}\psi_i'; \ddot M=\Sigma
_{ee}^{-1/2}\mathcal M(\Sigma_{ee}^{-1/2}\Psi)\Sigma_{ee}^{-1/2}$,
$\Gamma_p^h=(\gamma_{1p}^h, \gamma_{2p}^h,\ldots, \gamma_{N p}^h)'$.
\end{theorem}

%
\begin{remark} \label{remark-asy-beta}
In Appendix~B.3, we show
that the asymptotic expression of $\hat\beta-\beta$ in Theorem \ref{extenthmnew} can be expressed alternatively as
\begin{eqnarray*}
\hat\beta-\beta&=&\pmatrix{\operatorname{tr}
\bigl[\ddot MX_1\mathcal {M}(\widebar{\mathbb G})
X_1' \bigr]& \cdots& \operatorname{tr}\bigl[\ddot
MX_1\mathcal{M}(\widebar{\mathbb G}) X_K'
\bigr]
\vspace*{3pt}\cr
\vdots& \vdots& \vdots
\vspace*{3pt}\cr
\operatorname{tr}\bigl[\ddot MX_K\mathcal {M}(\widebar{\mathbb G})X_1' \bigr] & \cdots&\operatorname{tr}\bigl[\ddot
MX_K\mathcal{M}(\widebar{\mathbb G}) X_K'
\bigr]}^{-1}
\\[4pt]
&&{} \times\pmatrix{\operatorname{tr}\bigl[\ddot MX_1\mathcal{M}(\widebar{\mathbb G}) e' \bigr]
\vspace*{3pt}\cr
\vdots
\vspace*{3pt}\cr
\operatorname{tr}\bigl[\ddot MX _K\mathcal{M}(\widebar{\mathbb G}) e'\bigr]}
\\
&&{} +O_p
\bigl(T^{-3/2}\bigr)+O_p\bigl(N^{-1}T^{-1/2}
\bigr)+O_p\bigl(N^{-1/2}T^{-1}\bigr),
\end{eqnarray*}
where $X_k$ and $e$ are defined below (\ref{alternative1}) and
$\widebar{\mathbb G} =(1_T, \mathbb G)$. Notice $\ddot M$ is defined
as $\Sigma_{ee}^{-1/2}\mathcal{M}(\Sigma_{ee}^{-1/2}\Psi)\Sigma
_{ee}^{-1/2}$, which is equal to $\Sigma_{ee}^{-1/2}\mathcal
{M}(\Sigma
_{ee}^{-1/2}\Lambda)\Sigma_{ee}^{-1/2}$ since $\Lambda=(\Psi,
0_{N\times r_2})$ in the present context. In Appendix~B.3 of the
supplement~\cite{bai}, we also provide an intuitive explanation for this
alternative expression.
\end{remark}

Given Theorem \ref{extenthmnew} and Remark \ref{remark-asy-beta} we
have the following corollary.

%
\begin{corollary}[(Limiting distribution)]\label{extenthm2}
Under Assumptions \ref{assA}--\ref{assE}, if $\sqrt N/T\to0$, we have
\[
\sqrt{NT}(\hat\beta-\beta)\stackrel{d} {\rightarrow}N\bigl(0,{\widebar{
\mathcal P}}^{ -1}\bigr),
\]
where $\widebar{\mathcal P}=\lim_{N,T\to\infty} {\mathcal P}$,
and $\widebar{\mathcal P}$ is also the probability limit of
\[
\widebar{\mathcal P}=\mathop{\operatorname{plim}}_{N,T\to\infty
}
\frac{1} {NT}\pmatrix{\operatorname{tr}
\bigl[\ddot MX_1\mathcal {M}(\widebar{\mathbb G})
X_1' \bigr]& \cdots& \operatorname{tr}\bigl[\ddot
MX_1\mathcal{M}(\widebar{\mathbb G}) X_K'
\bigr]
\vspace*{3pt}\cr
\vdots& \vdots& \vdots
\vspace*{3pt}\cr
\operatorname{tr}\bigl[\ddot MX _K\mathcal {M}(\widebar{\mathbb G})X_1' \bigr] & \cdots&\operatorname{tr}\bigl[\ddot MX
_K\mathcal {M}(\widebar{\mathbb G}) X_K'
\bigr]}.
\]
\end{corollary}

%
\begin{remark} \label{remark1}
 Compared with the model in Section~\ref{sec2}, $\hat\beta$ is more efficient under the zero loading restrictions.
The reason is intuitive. In the previous model, only variations in
$v_{itx}$ provide information for $\beta$. But in the present case,
variations in $\gamma_{ik}^{h\prime}h_t$ of $x_{it}$ also provide
information for $\beta$. This can also be seen by comparing the
limiting variances of Corollaries \ref{comshkthm4} and \ref{extenthm2}.
Notice the projection matrix now only involves $\widebar{\mathbb G}$
instead of $\widebar{\mathbb F}$; and $\widebar{\mathbb G}$ is a
submatrix of $\widebar{\mathbb F}$. In addition, the covariance matrix
$\widebar{\mathcal P}$ can be estimated by the same method as in
estimating $\widebar\Omega$; see Remark \ref{remark2}.
\end{remark}

\section{Models with time-invariant regressors and common regressors}\label{sec4}
In this section, we extend the basic model in Section~\ref{sec2} to include
time-invariant regressors and common regressors. Examples of
time-invariant regressors include gender, race and education; and
examples for common regressors include price variables, unemployment
rate, or macroeconomic policy variables. These types of regressors are
important for empirical applications.

We first consider the model with only time-invariant regressors,
\begin{eqnarray}
\label{obseffeq2} y_{it}&=&\alpha_i+x_{it1}
\beta_1+x_{it2}\beta_2+\cdots +x_{itK}
\beta_K+\psi_i'g_t+
\phi_i'h_t+e_{it},
\nonumber
\\[-8pt]
\\[-8pt]
x_{itk}&=&\mu_{ik}+\gamma_{ik}^{g\prime}g_t+
\gamma_{ik}^{h\prime
}h_t+v_{itk}
\nonumber
\end{eqnarray}
for $k=1,2,\ldots,K$, where $g_t$ is an $r_1$-dimensional vector, and
$h_t$ is an \mbox{$r_2$-}dimen\-sional vector. Let $f_t=(g_t',h_t')'$, an
$r$-dimensional vector. The key point of model (\ref{obseffeq2}) is
that the $\phi_i$'s are known (but not zeros). We treat $\phi_i$ as new
added time-invariant regressors, whose coefficient $h_t$ is allowed to
be time-varying. The parameter of interest is still $\beta$. The
inference for $h_t$ is provided in Appendix~C.4 of the supplement~\cite{bai}.
The model in the previous section can be viewed as $\Phi=0$, where
$\Phi
=(\phi_1, \phi_2, \ldots, \phi_N)'$. However, the earlier
derivation is
not applicable here because now $\Phi$ is a general matrix with full
column rank, which provides more information (restrictions) on the
rotation matrix. Thus the number of restrictions required to eliminate
rotational indeterminacy is even fewer than in Section~\ref{sec3}. This point
can be seen in the next subsection.

We define the following notation for further analysis:
\begin{eqnarray*}
\Gamma_i^g&=&\bigl(\psi_i,
\gamma_{i1}^g,\ldots,\gamma_{iK}^g
\bigr),\qquad \Gamma_i^h=\bigl(\phi_i,
\gamma_{i1}^h,\ldots,\gamma_{iK}^h
\bigr),\qquad \Gamma_i=\bigl({\Gamma_i^g}',{
\Gamma_i^h}'\bigr)',
\\
\Phi&=&(\phi_1,\phi_2,\ldots,\phi_N)',
\qquad \Psi=(\psi_1,\psi_2,\ldots,\psi_N)',
\qquad \lambda_i=\bigl(\psi_i',
\phi_i'\bigr)',
\\
\Lambda&=&(\lambda_1,\lambda_2,\ldots,
\lambda_N)'.
\end{eqnarray*}
Then equation (\ref{obseffeq2}) has the same matrix expression as
(\ref
{comshkmaineq2}). Note that $\Lambda=[\Psi,\Phi]$ is the factor loading
matrix for the $N\times1$ vector $(y_{1t}, y_{2t},\ldots,y_{Nt})'$.

\subsection{Identification conditions}\label{sec4.1} We make the following
identification conditions, which we refer to as IO (\textit{Identification} conditions with partial \textit{Observable} fixed
effects), to emphasize the observed fixed effects:
\begin{longlist}[(IO1)]
\item[(IO1)] We partition the matrix $M_{f f}$ as
\[
M_{f f}=\lleft[
\matrix{M_{gg} &
M_{gh}
\vspace*{3pt}\cr
M_{hg} & M_{hh}}\rright]
\]
and impose $M_{gh}=0$ and $M_{gg}=I_{r_1}$;
\item[(IO2)] $\frac{1}{N}\Gamma^{g\prime}\Sigma_{\varepsilon
\varepsilon}^{-1}
\Gamma^g=D$, where $D$ is a diagonal matrix with its diagonal elements
distinct and arranged in descending order;
\item[(IO3)] $1_T'\mathbb G=0$ and $1_T'\mathbb H=0$.
\end{longlist}

In Appendix~C, we show that IO is sufficient for identification. These
restrictions can be imposed without loss of generality, as argued
formally in Appendix~C.3.
In addition, we make the following assumption.

\begin{ass}\label{assF} The loading matrix $\Lambda=[\Psi,\Phi]$ is of
full column rank.
\end{ass}

\subsection{Estimation}\label{sec4.2}
For clarity, in this subsection, we use $\Phi^*$ to denote the observed
value for $\Phi$. Recall that $\Sigma_{z z}=\Gamma M_{f f} \Gamma
'+\Sigma_{\varepsilon\varepsilon}$, where $\Gamma$ contains the factor
loading coefficients (including $\Phi$); $M_{f f}$ contains the
sub-blocks $M_{gg}$, $M_{gh}$ and $M_{hh}$; $\Sigma_{\varepsilon
\varepsilon}$ contains the heteroskedasticity coefficients. The
regression coefficient $\beta$ is contained in matrix $B$. The
maximization of the likelihood function is now subject to four sets of
restrictions, $M_{gh}=0$, $M_{g g}=I_{r_1}$, $\Phi=\Phi^*$ and
$\frac
{1}{N}\Gamma^{g\prime}\Sigma_{\varepsilon\varepsilon}^{-1}
\Gamma^g=D$. The likelihood function augmented with the Lagrange
multipliers is
\begin{eqnarray*}
\ln L&=&-\frac{1}{2N}\ln|\Sigma_{z z}|-\frac{1}{2N}
\operatorname {tr} \bigl[(I_N\otimes B)M_{z z}
\bigl(I_N\otimes B'\bigr)\Sigma_{z z}^{-1}
\bigr]+\operatorname {tr} [\Upsilon _1M_{gh} ]
\\
&&{} +\operatorname{tr} \bigl[\Upsilon_2(M_{g g}-I_{r_1})
\bigr]+\operatorname{tr} \biggl[\Upsilon_3 \biggl(\frac
{1}{N}
\Gamma^{g\prime}\Sigma_{\varepsilon\varepsilon}^{-1} \Gamma^g-D
\biggr) \biggr]+\operatorname{tr} \bigl[\Upsilon_4\bigl(\Phi-\Phi ^*
\bigr) \bigr],
\end{eqnarray*}
where $\Upsilon_1, \Upsilon_2, \Upsilon_3$ and $\Upsilon_4$ are all
Lagrange multipliers matrices; $\Upsilon_1$ is an $r_2\times r_1$
matrix; $\Upsilon_2$ is an $r_1\times r_1$ symmetric matrix; $\Upsilon
_3$ is an $r_1\times r_1$ symmetric matrix with all diagonal elements
zeros; $\Upsilon_4$ is an $r_2\times N$ matrix; and $\Sigma_{z
z}=\Gamma M_{f f}\Gamma'+\Sigma_{\varepsilon\varepsilon}$. Using the
same arguments in deriving $\Upsilon_1=0$ in Section~\ref{sec3}, we have
$\Upsilon_3=0$. Then the likelihood function is simplified as
\begin{eqnarray}
\label{obsobj} \ln L&=&-\frac{1}{2N}\ln|\Sigma_{z z}|-
\frac{1}{2N}\operatorname {tr} \bigl[(I_N\otimes
B)M_{z z}\bigl(I_N\otimes B'\bigr)
\Sigma_{z z}^{-1} \bigr]
\nonumber
\\[-8pt]
\\[-8pt]
&&{}+\operatorname{tr}[\Upsilon_1M_{g h}] +
\operatorname{tr}\bigl[\Upsilon_2(M_{g g}-I_{r_1})
\bigr]+\operatorname {tr}\bigl[\Upsilon_4\bigl(\Phi-\Phi^*\bigr)
\bigr].
\nonumber
\end{eqnarray}
The first order condition for $\Gamma$ gives
\[
\widehat M_{f f}\widehat\Gamma'\widehat\Sigma_{z z}^{-1}
\bigl[(I_N\otimes\widehat B)M_{z
z}\bigl(I_N\otimes
\widehat B{}'\bigr)-\widehat\Sigma_{z z}\bigr]\widehat\Sigma_{z z}^{-1}=W',
\]
where $W$ is defined in (\ref{extenmaineq1}). Pre-multiplying $\widehat M_{f f}^{-1}$ and post-multiplying $\widehat\Gamma$, and by the special
structures of $W$ and $\widehat\Gamma$, we have
\begin{eqnarray*}
&& \frac{1}{N}\widehat\Gamma'\widehat\Sigma_{z z}^{-1}
\bigl[(I_N\otimes\widehat B)M_{z
z}\bigl(I_N\otimes
\widehat B{}'\bigr)-\widehat\Sigma_{z z}\bigr]\widehat\Sigma_{z z}^{-1}\widehat\Gamma
\\
&&\qquad = -
\lleft[
\matrix{0_{r_1\times r_1} & 0_{r_1\times r_2}
\vspace*{3pt}\cr
\dfrac{1} N \widehat M_{hh}^{-1}\Upsilon_4'
\widehat\Psi& \dfrac{1} N\widehat M_{hh}^{-1}\Upsilon
_4'\Phi}\rright].
\end{eqnarray*}
But the first order condition for $M_{f f}$ gives
\begin{equation}
\label{obseffeq0715}\quad \frac{1}{N}\widehat\Gamma'\widehat\Sigma
_{z z}^{-1}\bigl[(I_N\otimes\widehat B)M_{z z}\bigl(I_N\otimes\widehat B{}'\bigr)-\widehat\Sigma _{z z}\bigr]\widehat\Sigma_{z z}^{-1}\widehat\Gamma=
\lleft[\matrix{
\Upsilon_2 &
\Upsilon_1'
\vspace*{3pt}\cr
\Upsilon_1 & 0_{r_2\times r_2}}
\rright].
\end{equation}
Comparing\vspace*{1pt} the proceeding two results and noting that the left-hand side
is a
symmetric matrix,\vspace*{1pt} we have $\widehat\Gamma'\widehat\Sigma_{z
z}^{-1}[(I_N\otimes\widehat B)M_{z z}(I_N\otimes\widehat B{}')-\widehat\Sigma
_{z
z}]\widehat\Sigma_{z z}^{-1}\widehat\Gamma=0$. But $\widehat\Gamma'\widehat\Sigma_{z
z}^{-1}$ can be replaced by $\widehat\Gamma'\widehat\Sigma_{\varepsilon
\varepsilon}^{-1}$; see (S.2) in the Appendix. Thus
\begin{equation}
\label{obseffeq3}\widehat\Gamma'\widehat\Sigma _{\varepsilon
\varepsilon}^{-1}
\bigl[(I_N\otimes\widehat B)M_{z z}\bigl(I_N\otimes
\widehat B{}'\bigr)-\widehat\Sigma_{z z}\bigr]\widehat\Sigma_{\varepsilon\varepsilon}^{-1}\widehat\Gamma=0.
\end{equation}
The above result implies that $\Upsilon_1=0$, $\Upsilon_2=0$,
$\Upsilon
_4'\widehat\Psi=0$ and $\Upsilon_4'\Phi=0$.

The first order condition for $\Sigma_{ii}$ is the same as (\ref
{extenmaineq5}), that is,
\begin{eqnarray}
\label{obseffeq4} && \widehat B M_{z z}^{j j}\widehat B{}'-\widehat\Sigma_{z z}^{j
j}-\widehat\Gamma_j'\widehat G\sum_{i=1}^N \widehat\Gamma_i
\widehat\Sigma _{ii}^{-1} \bigl(\widehat B M_{z z}^{i j}
\widehat B{}'-\widehat\Sigma_{z z}^{i j} \bigr)
\nonumber
\\[-8pt]
\\[-8pt]
&&\qquad{}-\sum_{i=1}^N \bigl(\widehat BM_{z z}^{ji}\widehat B{}'-\widehat\Sigma_{z
z}^{ji}
\bigr)\widehat\Sigma_{ii}^{-1}\widehat\Gamma_i'
\widehat G\widehat\Gamma_j=\mathbb {W},
\nonumber
\end{eqnarray}
where $\mathbb{W}$ is defined following (\ref{comshkmaineq7}).

The first order condition on $\beta$ is the same as (\ref
{extenmainbeta}), that is,
\begin{equation}
\label{obseffeq1} \frac{1}{NT}\sum_{i=1}^N
\sum_{t=1}^T\widehat\Sigma_{iie}^{-1}
\Biggl\{(\dot y_{it}- \dot x_{it}\hat\beta) -\hat
\lambda_i'\widehat G\sum_{j=1}^N
\widehat\Gamma_j\widehat\Sigma _{jj}^{-1}
\lleft[
\matrix{\dot y_{jt} -\dot x_{jt}\hat
\beta
\vspace*{3pt}\cr
\dot x_{jt}'}\rright]
\Biggr\}\dot x_{it}=0.\hspace*{-35pt}
\end{equation}

We need an additional identify for the theoretical analysis in the
Appendix. The preceding analysis shows that $\frac{1}{N}\Upsilon
_4'\widehat\Psi=0$ and $\frac{1}{N}\Upsilon_4'\Phi=0$. They imply
\begin{equation}
\label{obseffeq5} \frac{1}{N}\sum_{i=1}^N
\sum_{j=1}^N\widehat G_{2}\widehat\Gamma_i\widehat\Sigma _{ii}^{-1}\bigl(\widehat BM_{z z}^{i j}\widehat B{}'-\widehat\Sigma_{z z}^{i j}
\bigr)\widehat\Sigma_{j j}^{-1}I _{K+1}^{ 1}
\hat\lambda_j'=0,
\end{equation}
where $\hat\lambda_j=(\hat\psi_i',\phi_i')'$.

\subsection{Asymptotic properties}\label{sec4.3} The asymptotic representation for
$\hat\beta-\beta$ is:
%
\begin{proposition}\label{obseffthm1}
Under Assumptions \ref{assA}--\ref{assD} and \ref{assF}, and under the identification condition
\textup{IO}, we have
\begin{eqnarray*}
\mathcal{Q}^0(\hat\beta-\beta)&=& \frac{1}{NT}\sum
_{i=1}^N\sum_{t=1}^T
\Sigma_{iie}^{-1}e_{it}v_{itx}+
\frac{1}{NT} \sum_{i=1}^N\sum
_{t=1}^T\Sigma_{iie}^{-1}
\gamma_{ix}^{h\prime} h_te_{it}
\\
&&{}-\frac{1}{NT}\sum_{i=1}^N\sum
_{t=1}^T\Sigma_{iie}^{-1}
\lambda _i'\Pi _{\lambda\lambda}^{-1} \Biggl(
\frac{1}{N} \sum_{j=1}^N
\lambda_j'\Sigma_{j je}^{-1}
\gamma_{jx}^{h\prime} \Biggr) h_te_{it}
\\
&&{}+O_p\bigl(T^{-3/2}\bigr)+O_p
\bigl(N^{-1}T^{-1/2}\bigr)+O_p\bigl(N^{-1/2}T^{-1}
\bigr),
\end{eqnarray*}
where $\mathcal{Q}^0$ is a $K\times K$ symmetric matrix with its
$(p,q)$ element equal to
$\frac{1} N\operatorname{tr}[M_{h h}\Gamma_p^{h\prime}\ddot M\Gamma
_q^h]+\frac{1} N\sum_{i=1}^N\Sigma_{iie}^{-1}\Sigma_{iix}^{(p,q)}$; $\ddot M=\Sigma
_{ee}^{-1/2}\mathcal M(\Sigma_{ee}^{-1/2}\Lambda)\Sigma_{ee}^{-1/2}$;
$\Gamma_p^h=[\gamma_{1p}^h, \gamma_{2p}^h,\ldots, \gamma
_{Np}^h]'$; $\Pi
_{\lambda\lambda}=\frac{1} N\sum_{i=1}^N\lambda_i\Sigma
_{iie}^{-1}\lambda
_i'$; and $\gamma_{jx}^h=[\gamma_{j1}^h, \gamma_{j2}^h,\ldots,
\gamma_{jK}^h]$.
\end{proposition}

Proposition \ref{obseffthm1} is derived under the identification
conditions IO. In Appendix~C.3, we show that for any set of factors and
factor loadings $(\psi_i, \gamma_{ik},\break  g_t, h_t)$, we can always\vadjust{\goodbreak}
transform it to another set $(\psi_i^\star, \gamma_{ik}^\star,
g_t^\star, h_t^\star)$ which satisfies IO, and at the same time, still
maintains
the observability of $\Phi$ (i.e., $\Phi$ is untransformed). This is in
agreement with the Lagrange multiplier analysis, in which $\Upsilon
_j=0$ ($j=1,2,3)$, but the multiplier for $\Phi=\Phi^*$ is nonzero.
Using the relationship between the two sets, we can generalize
Proposition \ref{obseffthm1} into the following theorem, which does not
depend on IO.

%
\begin{theorem}\label{obseffthmnew}
Under Assumptions \ref{assA}--\ref{assD} and \ref{assF}, we have
\begin{eqnarray*}
\mathcal{Q}(\hat\beta-\beta)&=& \frac{1}{NT}\sum
_{i=1}^N\sum_{t=1}^T
\Sigma_{iie}^{-1} e_{it}v_{itx}+
\frac{1}{NT}\sum_{i=1}^N\sum
_{t=1}^T\Sigma _{iie}^{-1}
\gamma_{ix}^{h\prime} h_t^\star
e_{it}
\\
&&{}-\frac{1}{NT}\sum_{i=1}^N\sum
_{t=1}^T\Sigma_{iie}^{-1}
\lambda _i'\Pi _{\lambda\lambda}^{-1} \Biggl(
\frac{1}{N}\sum_{j=1}^N
\lambda_j' \Sigma_{j je}^{-1}
\gamma_{jx}^{h\prime} \Biggr)h_t^\star
e_{it}
\\
&&{}+O_p\bigl(T^{-3/2}\bigr)+O_p
\bigl(N^{-1}T^{-1/2}\bigr)+O_p\bigl(N^{-1/2}T^{-1}
\bigr),
\end{eqnarray*}
where
\[
h_t^\star=\dot h_t-\dot{\mathbb
H}'\dot{\mathbb G}\bigl(\dot{\mathbb G}'\dot {\mathbb
G}\bigr)^{-1}\dot g_t;
\]
$\mathcal{Q}$ is a $K\times K$ symmetric matrix with its $(p,q)$
element equal to
\[
\frac{1} {NT}\operatorname{tr}\bigl[\ddot M\Gamma_q^h
\mathbb H'\mathcal M(\widebar{\mathbb G})\mathbb H
\Gamma_p^{h\prime}\bigr]+\frac{1} N\sum
_{i=1}^N\Sigma _{iie}^{-1}
\Sigma_{iix}^{(p,q)}
\]
and $\ddot M$, $\Gamma_p^h$ and $\Pi_{\lambda\lambda}$ are defined in
Proposition \ref{obseffthm1}.
\end{theorem}

%
\begin{remark}\label{remark3}
In Appendix~C.3 we show that the asymptotic expression of $\hat\beta
-\beta$ in Theorem \ref{obseffthmnew} can be expressed alternatively as
\begin{eqnarray*}
\hat\beta-\beta&=&\pmatrix{\operatorname{tr}
\bigl[\ddot MX_1\mathcal {M}(\widebar{\mathbb G})
X_1' \bigr]& \cdots& \operatorname{tr}\bigl[\ddot
MX_1\mathcal{M}(\widebar{\mathbb G}) X_K'
\bigr]
\vspace*{3pt}\cr
\vdots& \vdots& \vdots
\vspace*{3pt}\cr
\operatorname{tr}\bigl[\ddot MX _K\mathcal {M}(\widebar{\mathbb G})X_1' \bigr] & \cdots&\operatorname{tr}\bigl[\ddot MX
_K\mathcal {M}(\widebar{\mathbb G}) X _K'
\bigr]}^{-1}
\\[4pt]
&&{}\times\pmatrix{\operatorname{tr}\bigl[\ddot
MX_1\mathcal{M}(\widebar{\mathbb G}) e' \bigr]
\vspace*{3pt}\cr
\vdots
\vspace*{3pt}\cr
\operatorname{tr}\bigl[\ddot MX _K\mathcal{M}(\widebar{\mathbb G}) e'\bigr]}+O_p
\bigl(T^{-3/2}\bigr)
\\
&&{}+O_p\bigl(N^{-1}T^{-1/2}
\bigr)+O_p\bigl(N^{-1/2}T^{-1}\bigr),
\end{eqnarray*}
where $X_k$ and $e$ are defined below (\ref{alternative1}) and
$\widebar{\mathbb G}=(1_T, \mathbb G)$. We also show in Appendix~C.3
that this alternative expression has an intuitive explanation.
\end{remark}

From Theorem \ref{obseffthmnew}, we obtain the following corollary.

%
\begin{corollary}\label{obseffcor}
Under the conditions of Theorem \ref{obseffthm1}, if $\sqrt N/T\to0$,
we have
\[
\sqrt{NT}(\hat\beta-\beta)\stackrel{d} {\rightarrow}N\bigl(0,{\widebar{\mathcal Q}}^{ -1}\bigr),
\]
where $\widebar{\mathcal Q}=\lim_{N,T\to\infty}\mathcal Q$,
which has an alternative expression
\[
\widebar{\mathcal Q}=\mathop{\operatorname{plim}}_{N,T\to\infty
}
\frac{1} {NT}\pmatrix{\operatorname{tr}
\bigl[\ddot MX_1\mathcal {M}(\widebar{\mathbb G})
X_1' \bigr]& \cdots& \operatorname{tr}\bigl[\ddot
MX_1\mathcal{M}(\widebar{\mathbb G}) X_K'
\bigr]
\vspace*{3pt}\cr
\vdots& \vdots& \vdots
\vspace*{3pt}\cr
\operatorname{tr}\bigl[\ddot MX _K\mathcal {M}(\widebar{\mathbb G})X_1' \bigr] & \cdots&\operatorname{tr}\bigl[\ddot
MX_K\mathcal {M}(\widebar{\mathbb G}) X _K'
\bigr]}.
\]
\end{corollary}

%
\begin{remark}
Compared with the model in Section~\ref{sec2}, $\hat\beta$ is more efficient
with observable fixed effects (time-invariant regressors). The reason
is provided in Remark \ref{remark1}.
\end{remark}

\subsection{Models with time-invariant regressors and common regressors}\label{sec4.4}
In this subsection, we consider the joint presence of time-invariant
regressors and common regressors. Consider the following model:
\begin{eqnarray}
\label{obseffeq7} y_{it}&=&x_{it1}\beta_1+x_{it2}
\beta_2+\cdots+x_{it K}\beta _K+
\psi_i'g_t+\phi_i'h_t+
\kappa_i'd_t+e_{it},
\nonumber
\\[-8pt]
\\[-8pt]
x_{itk}&=&\gamma_{ik}^{g\prime}g_t+
\gamma_{ik}^{h\prime}h_t+\gamma _{ik}^{d\prime}\,d_t+v_{itk}
\nonumber
\end{eqnarray}
for $k=1,2,\ldots,K$, where $g_t$, $h_t$ and $d_t$ are $r_1\times1$,
$r_2\times1$ and $r_3\times1$ vectors, respectively.
A key feature of model (\ref{obseffeq7}) is that $d_t$ and $\phi_i$
are observable for all $i$ and $t$. We call $\phi_i$ the time-invariant
regressors because they are invariant over time and $d_t$ the common
regressors because they are the same for all the cross-sectional units.
In this model, the time-invariant regressors have time-varying
coefficients, and the common regressors have heterogeneous
(individual-dependent) coefficients. If $d_t\equiv1$, $\kappa_i$ plays
the role of $\alpha_i$ in (\ref{obseffeq2}). So the model here is
more general.

Similar to the previous subsection, we make the following assumption:

\begin{ass}\label{assG} The matrices $(\Psi, \Phi, \mathrm{K})$ and
$(\mathbb G, \mathbb H, \mathbb D)$ are both of full column rank, where
$\mathrm{K}=(\kappa_1, \kappa_2,\ldots, \kappa_N)'$ and $\mathbb D=(d_1,
d_2,\ldots, d_T)'$.
\end{ass}

Let $\lambda_i=(\psi_i',\phi_i')'$, $\gamma_{ik}=(\gamma
_{ik}^{g\prime
},\gamma_{ik}^{h\prime})'$ and $\delta_i=(\kappa_i,\gamma
_{ik}^d)$. The
model can be written as
\[
\lleft[\matrix{1 & -\beta'
\vspace*{3pt}\cr
0 & I_K}\rright]z_{it}=
\Gamma_i'f_t +\delta_i'd_t+
\varepsilon_{it},
\]
where $z_{it}, \Gamma_i, \varepsilon_{it}$ are defined in Section~\ref{sec2};
Let $\Delta=(\delta_1,\delta_2,\ldots, \delta_N)'$. Then
\begin{equation}
\label{obseffeq8}(I_N\otimes B)z_t-\Delta d_t=
\Gamma f_t+\varepsilon_t,
\end{equation}
where the symbols $\Gamma, z_t, B, \varepsilon_t$ are defiend in
Section~\ref{sec2}.

The likelihood function can be written as
\[
\ln L =-\frac{1} {2N}\ln|\Sigma_{z z}|-\frac{1} {2NT}\sum
_{t=1}^T\bigl[(I_N\otimes
B)z_t-\Delta d_t\bigr]'
\Sigma_{z z}^{-1}\bigl[(I_N\otimes
B)z_t-\Delta d_t\bigr].
\]
Take $\Sigma_{z z}$ and $\beta$ as given. $\Delta$ maximizes the above
function at
\[
\widehat\Delta=(I_N\otimes B) \Biggl(\sum_{s=1}^Tz_sd_s'
\Biggr) \Biggl(\sum_{s=1}^T\,d_sd_s'
\Biggr)^{-1}.
\]
Substituting $\widehat\Delta$ into the above likelihood function, we obtain
the concentrated likelihood function
\[
\ln L=-\frac{1} {2N}\ln|\Sigma_{z z}|-\frac{1} {2NT}
\operatorname {tr} \bigl[(I_N\otimes B)Z\mathcal{M}(\mathbb
D)Z'\bigl(I_N\otimes B'\bigr)
\Sigma_{z z}^{-1} \bigr],
\]
where $Z=(z_1,z_2,\ldots, z_T)$, $\mathbb{D}=(d_1,d_2,\ldots,d_T)'$ and
$\mathcal{M}(\mathbb D)=I_T-\mathbb D({\mathbb D}'\mathbb
D)^{-1}{\mathbb D}'$, a projection matrix. Consider (\ref{obseffeq8}),
which is equivalent to
\[
(I_N\otimes B)Z=\Gamma\mathbb F'+\Delta\mathbb
D' +\varepsilon,
\]
where $\varepsilon=(\varepsilon_1, \varepsilon_2,\ldots,
\varepsilon
_T)$. Post-multiplying $\mathcal M(\mathbb D)$ on both sides, we have
\[
(I_N\otimes B)Z\mathcal M(\mathbb D)=\Gamma\mathbb F'
\mathcal M(\mathbb D) +\varepsilon\mathcal M(\mathbb D).
\]
If we treat $Z\mathcal M(\mathbb D)$ as the new observable data,
$\mathbb F'\mathcal M(\mathbb D)$ as the new unobservable factors, the
preceding equation can be viewed as a special case of (\ref
{obseffeq2}). Invoking Theorem \ref{obseffthmnew}, which does not need
IO [the factors $\mathbb F'\mathcal M(\mathbb D)$ may not satisfy IO],
we have the following theorem:

%
\begin{theorem}\label{obseffthm3}
Under Assumptions \ref{assA}--\ref{assD} and \ref{assG}, the asymptotic representation of $\hat
\beta$ in the presence of time invariant and common regressors is
\begin{eqnarray*}
\mathcal{R}(\hat\beta-\beta)&=& \frac{1}{NT}\sum
_{i=1}^N\sum_{t=1}^T
\Sigma_{iie}^{-1}e_{it}v_{itx}+
\frac{1}{NT}\sum_{i=1}^N \sum
_{t=1}^T\Sigma_{iie}^{-1}
\gamma_{ix}^{h\prime}h_t^\star
e_{it}
\\
&&{}-\frac{1}{NT}\sum_{i=1}^N\sum
_{t=1}^T\Sigma_{iie}^{-1}
\lambda _i'\Pi _{\lambda\lambda}^{-1}
\frac{1}{N}\sum_{j=1}^N
\lambda_j' \Sigma_{j je}^{-1}
\gamma_{jx}^{h\prime}h_t^\star
e_{it}
\\
&&{}+O_p\bigl(T^{-3/2}\bigr)+O_p
\bigl(N^{-1}T^{-1/2}\bigr)+O_p\bigl(N^{-1/2}T^{-1}
\bigr),
\end{eqnarray*}
where
\[
h_t^\star=h_t-\mathbb H'\mathbb
D\bigl(\mathbb D'\mathbb D\bigr)^{-1}\,d_t-
\mathbb H'\mathcal M(\mathbb D)\mathbb G\bigl[\mathbb
G'\mathcal M(\mathbb D)\mathbb G\bigr]^{-1}
\bigl(g_t-\mathbb G'\mathbb D\bigl(\mathbb
D'\mathbb D\bigr)^{-1}\,d_t\bigr);
\]
$\mathcal{R}$ is a $K\times K$ symmetric matrix with its $(p,q)$
element equal to
\[
\frac{1} {NT}\operatorname{tr}\bigl[\ddot M\Gamma_q^h
\mathbb H'\mathcal M(\mathbb B)\mathbb H\Gamma_p^{h\prime}
\bigr]+\frac{1} N\sum_{i=1}^N
\Sigma_{iie}^{-1}\Sigma _{iix}^{(p,q)},
\]
where $b_t=(g_t', d_t')'$ and $\mathbb B=(b_1, b_2, \ldots,
b_T)'=(\mathbb G, \mathbb D)$, a matrix of $T\times(r_1+r_3)$
dimension; $\ddot M=\Sigma_{ee}^{-1/2}\mathcal M(\Sigma
_{ee}^{-1/2}\Lambda)\Sigma_{ee}^{-1/2}; \Gamma_p^h=(\gamma_{1p}^h,
\gamma_{2p}^h,\ldots, \gamma_{N p}^h)'$; $\Pi_{\lambda\lambda
}=\break \frac{1} N
\sum_{i=1}^N\lambda_i\Sigma_{iie}^{-1}\lambda_i'$.
\end{theorem}

%
\begin{remark} \label{remark-invariant-common}
The asymptotic expression of $\hat\beta-\beta$ can be alternatively
expressed as
\begin{eqnarray*}
\hat\beta-\beta&=&\pmatrix{\operatorname{tr}
\bigl[\ddot MX_1\mathcal{M}(\mathbb B) X_1'
\bigr] & \cdots& \operatorname{tr}\bigl[\ddot MX_1\mathcal{M}(\mathbb
B) X_K' \bigr]
\vspace*{3pt}\cr
\vdots& \vdots& \vdots
\vspace*{3pt}\cr
\operatorname{tr}\bigl[\ddot MX_K\mathcal {M}(\mathbb
B)X_1' \bigr] & \cdots& \operatorname{tr}\bigl[\ddot
MX_K\mathcal{M}(\mathbb B) X_K' \bigr]}^{-1}
\\[4pt]
&&{} \times\pmatrix{\operatorname{tr}\bigl[\ddot
MX_1\mathcal{M}(\mathbb B) e' \bigr]
\vspace*{3pt}\cr
\vdots
\vspace*{3pt}\cr
\operatorname{tr}\bigl[\ddot MX_K\mathcal{M}(\mathbb B)
e'\bigr]}
\\
&&{} +O_p
\bigl(T^{-3/2}\bigr)+O_p\bigl(N^{-1}T^{-1/2} \bigr)+O_p\bigl(N^{-1/2}T^{-1}\bigr).
\end{eqnarray*}
If $\mathbb D=1_T$, the above asymptotic result reduces to the one in
Theorem \ref{obseffthmnew} since $\mathbb B=(1_T, \mathbb G)=\widebar{\mathbb G}$.
\end{remark}

Given Theorem \ref{obseffthm3} and Remark \ref
{remark-invariant-common}, we have the following corollary.

%
\begin{corollary}\label{corsub44}
Under Assumptions \ref{assA}--\ref{assD} and \ref{assG}, if $\sqrt N/T \to0$, then
\[
\sqrt{NT}(\hat\beta-\beta)\stackrel{d} {\rightarrow}N\bigl(0,{\widebar{\mathcal {R}}}^{ -1}\bigr),
\]
where $\widebar{\mathcal R}=\lim_{N,T\to\infty} \mathcal R$,
and $\widebar{\mathcal R}$ can also be expressed as
\[
\widebar{\mathcal R}=\mathop{\operatorname{plim}}_{N,T\to\infty
}
\frac{1} {NT}\pmatrix{\operatorname{tr}
\bigl[\ddot MX_1\mathcal{M}(\mathbb B) X_1'
\bigr] & \cdots& \operatorname{tr}\bigl[\ddot MX_1\mathcal{M}(\mathbb
B) X_K' \bigr]
\vspace*{3pt}\cr
\vdots& \vdots& \vdots
\vspace*{3pt}\cr
\operatorname{tr}\bigl[\ddot MX_K\mathcal {M}(\mathbb
B)X_1' \bigr] & \cdots&\operatorname{tr}\bigl[\ddot
MX_K\mathcal{M}(\mathbb B) X_K' \bigr]}.
\]
\end{corollary}

\section{Computing algorithm}\label{sec5}
To estimate the model by the maximum likelihood method, we adapt the
ECM (expectation and conditional maximization) procedures of \cite
{meng1993maximum}. More specifically, in the M-step we split the
parameter $\theta=(\beta, \Gamma, \Sigma_{\varepsilon\varepsilon},
M_{f f})$ into two blocks, $\theta_1=(\Gamma, \Sigma_{\varepsilon
\varepsilon}, M_{f f})$ and $\theta_2=\beta$, and update $\theta
_1^{(k)}$ to $\theta_1^{(k+1)}$ given $\theta_2^{(k)}$ and then update
$\theta_2^{(k)}$ to $\theta_2^{(k+1)}$ given $\theta_1^{(k+1)}$, where
$\theta^{(k)}$ is the estimated value at the $k$th iteration. In this
section, we only state the iterating formulas for basic models. The
iterating formulas for the models in Sections~\ref{sec3} and \ref{sec4} can be found in
Appendix~E of \cite{bai}. In Appendix~E, we also show that the iterated
EM solutions satisfy the first order conditions. So the EM estimators
are at least locally optimal.

In the basic model, $M_{f f}=I_r$. So the parameters to be estimated
reduce to $\theta=(\beta, \Gamma, \Sigma_{\varepsilon\varepsilon})$.
Let $\theta^{(k)}=(\beta^{(k)}, \Gamma^{(k)}, \Sigma_{\varepsilon
\varepsilon}^{(k)})$ be the estimated value at the $k$th iteration. We
update $\Gamma^{(k)}$ according to
\begin{equation}
\label{51}\Gamma^{(k+1)}= \Biggl[\frac{1} T\sum
_{t=1}^TE\bigl(z_tf_t'|Z,
\theta^{(k)}\bigr) \Biggr] \Biggl[\frac{1} T\sum
_{t=1}^TE\bigl(f_tf_t'|Z,
\theta^{(k)}\bigr) \Biggr]^{-1},
\end{equation}
where
\begin{eqnarray}
&& \frac{1} T\sum_{t=1}^TE
\bigl(f_t f_t'|Z,\theta^{(k)}
\bigr)\nonumber
\\[-2pt]
&&\qquad =I_r-\Gamma ^{(k)\prime
}\bigl(\Sigma_{z z}^{(k)}
\bigr)^{-1}\Gamma^{(k)}\label{52}
\\
&&\quad\qquad{} +\Gamma^{(k)\prime}\bigl(\Sigma_{z z}^{(k)}
\bigr)^{-1}\bigl(I_N\otimes B^{(k)}
\bigr)M_{z
z}\bigl(I_N\otimes B^{(k)\prime}\bigr) \bigl(
\Sigma_{z z}^{(k)}\bigr)^{-1}\Gamma ^{(k)},
\nonumber
\\[-2pt]
\label{53} &&\frac{1} T\sum
_{t=1}^TE\bigl(z_tf_t'|Z,
\theta ^{(k)}\bigr)=\bigl(I_N\otimes B^{(k)}
\bigr)M_{z z}\bigl(I_N\otimes B^{(k)\prime
}\bigr) \bigl(
\Sigma _{z z}^{(k)}\bigr)^{-1}\Gamma^{(k)}
\end{eqnarray}
with $\Sigma_{z z}^{(k)}=\Gamma^{(k)}\Gamma^{(k)\prime}+\Sigma
_{\varepsilon\varepsilon}^{(k)}$. We update $\Sigma_{\varepsilon
\varepsilon}^{(k)}$ and $\beta^{(k)}$ according to
\begin{eqnarray}
\Sigma_{\varepsilon\varepsilon}^{(k+1)}&=&\operatorname{Dg} \bigl\{
\bigl(I_{N(K+1)}-\Gamma^{(k+1)}\Gamma^{(k)\prime}\bigl(
\Sigma_{z
z}^{(k)}\bigr)^{-1} \bigr)
\nonumber
\\[-9pt]\label{54}
\\[-9pt]
&&\hspace*{27pt}{} \times\bigl(I_N\otimes B^{(k)}\bigr)M_{z z}
\bigl(I_N\otimes B^{(k)\prime}\bigr) \bigr\},
\nonumber
\\
\beta^{(k+1)}&=& \Biggl(\sum_{i=1}^N
\sum_{t=1}^T\dot x_{it}'
\bigl(\Sigma _{iie}^{(k+1)}\bigr)^{-1}\dot
x_{it} \Biggr)^{-1}
\nonumber
\\[-9pt]\label{55}
\\[-9pt]
&&{}\times \Biggl(\sum_{i=1}^N\sum
_{t=1}^T\dot x_{it}'
\bigl(\Sigma _{iie}^{(k+1)}\bigr)^{-1}\bigl(\dot
y_{it}-\lambda_i^{(k+1)\prime
}f_t^{(k)}
\bigr) \Biggr),
\nonumber
\end{eqnarray}
where $f_t^{(k)}$ is the transpose of the $t$th row of
\[
\mathbb F^{(k)}=E\bigl(\mathbb F|Z,\theta^{(k)}\bigr)=\dot
Z'\bigl(I_N\otimes B^{(k)\prime}\bigr) \bigl(
\Sigma_{z z}^{(k)}\bigr)^{-1}\Gamma^{(k)},
\]
where $\dot Z=(\dot z_1,\dot z_2,\ldots,\dot z_T)$ with $\dot
z_t=z_t-\frac{1} T\sum_{s=1}^Tz_s$; $\operatorname{Dg}(\cdot)$ is the operator
that sets the entries of its argument to zeros if the counterparts of
$E(\varepsilon_t\varepsilon_t')$ are zeros.

Putting together, we obtain $\theta^{(k+1)}=(\Gamma^{(k+1)}, \beta
^{(k+1)}, \Sigma_{\varepsilon\varepsilon}^{(k+1)})$. The above
iteration continues until $\|\theta^{(k+1)}-\theta^{(k)}\|$ is smaller
than a preset error tolerance. The initial values use the iterated PC
estimators of \cite{bai2009panel}.

\section{Finite sample properties}\label{sec6}

In this section, we consider the finite sample properties of the MLE.
Data are generated according to
\begin{eqnarray}
\label{datagen} y_{it}&=&\alpha_i+x_{it1}\beta
_1+x_{it2}\beta_2+\psi_i'g_t+
\phi_i'h_t+\kappa_i'd_t+e_{it},
\nonumber
\\[-8pt]
\\[-8pt]
x_{itk}&=&\mu_{ik}+\gamma_{ik}^{g\prime}g_t+
\gamma_{ik}^{h\prime
}h_t+\gamma_{ik}^{d\prime}\,d_t+v_{itk},
\qquad k=1,2.
\nonumber
\end{eqnarray}
The dimensions of $g_t, h_t, d_t$ are each fixed to 1. We set $\beta
_1=1$ and $\beta_2=2$. We consider four types of DGP (data generating
process), which correspond to the four models considered in the paper.
\begin{longlist}[DGP1:]
\item[DGP1:] $\phi_i, \kappa_i, \gamma_{ik}^h$ and
$\gamma
_{ik}^d$ are fixed to zeros; $\alpha_i, \mu_{ik}, \psi_i$ and $g_t$ are
generated from $N(0,1)$ and $\gamma_{ik}^g=\psi_i+N(0,1)$.

\item[DGP2:] $\phi_i, \kappa_i$ and $\gamma_{ik}^d$ are fixed
to zeros; $\alpha_i, \mu_{ik}, \psi_i$, $\gamma_{ik}^h, g_t$ and $h_t$
are generated from $N(0,1)$; $\gamma_{ik}^g=\psi_i+N(0,1)$.

\item[DGP3:] $\kappa_i$ and $\gamma_{ik}^d$ are fixed to
zeros; $\alpha_i, \mu_{ik}, \psi_i, \phi_i, g_t$ and $h_t$ are
generated from $N(0,1)$; $\gamma_{ik}^g=\psi_i+N(0,1)$ and $\gamma
_{ik}^h=\phi_i+N(0,1)$. Here $\phi_i$ is observable.

\item[DGP4:] $\alpha_i, \mu_{ik}, \psi_i, \phi_i, \kappa_i,
g_t$ and $h_t$ are generated from $N(0,1)$; $d_t=1+N(0,1)$, $\gamma
_{ik}^g=\psi_i+N(0,1)$, $\gamma_{ik}^h=\phi_i+N(0,1)$ and $\gamma
_{ik}^d=\kappa_i+N(0,1)$. Here $\phi_i$ and $d_t$ are observable.
\end{longlist}
Using the method of writing (\ref{comshkmaineq2}), we can rewrite
(\ref
{datagen}) as
\[
(I_N\otimes B)z_t=\mu+ L\varsigma_t+
\varepsilon_t,
\]
where $\varsigma_t=g_t$ for DGP1; $\varsigma_t=(g_t, h_t)'$ for DGP2
and DGP3; $\varsigma_t=(g_t, h_t, d_t)'$ for DGP4, and $L$ is the
corresponding loadings matrix. Let $\iota_i'$ be the $i$th row of $L$.
We generate the cross-sectional heteroscedasticity $\Xi$, an
$N(K+1)\times1$ vector, according to $\Xi_i=\frac{\eta_i} {1-\eta_i}
\iota_i'\iota_i, i=1, 2,\ldots, N(K+1)$,
where $\eta_i$ is drawn from $U[u, 1-u]$ with $u=0.1$. A similar way of
generating heteroscedasticity is also used in \cite{breitung2011gls}
and \cite{doz2006quasi}. Let $\Upsilon=\operatorname{diag}(\Upsilon
_1, \Upsilon_2,\ldots, \Upsilon_N)$ be an $N(K+1)\times N(K+1)$ block diagonal matrix,
in which $\Upsilon_i=\operatorname{diag}\{1$, $ (M_i'M_i)^{-1/2}M_i\}
$ with $M_i$
being a $K\times K$ standard normal random matrix for each $i$. Once
$\Upsilon$ is generated, the error term~$\varepsilon_t$, which is
defined as $(\varepsilon_{1t}',\varepsilon_{2t}',\ldots, \varepsilon
_{Nt}')'$ with $\varepsilon_{it}=(e_{it},v_{it1},v_{it2})'$, is
calculated by $\varepsilon_t=\sqrt{\operatorname{diag}(\Xi
)}\Upsilon\epsilon_t$,
where $\epsilon_t$ is an $N(K+1)\times1$ vector with all its elements
being i.i.d. $(\chi_2^2-2)/2$, where $\chi_2^2$ denotes the chi-squared
distribution with two freedom degrees, which is normalized to mean zero
and variance one. Additional simulation results for normal and
student-$t$ errors are given in Appendix~D. Once $\varepsilon_t$ is
obtained, we use
\[
z_t=(I_N\otimes B)^{-1}(\mu+L
\varsigma_t+\varepsilon_t)
\]
to yield the observable data.

In the basic model, the number of factors is determined by
\begin{equation}
\label{eq0} \hat r=\mathop{\operatorname{argmin}}_{0\le m\le r_{\max}} \operatorname{IC}(m)
\end{equation}
with
\[
\operatorname{IC}(m)=\frac{1} {N\widebar K}\ln \bigl|\widehat\Gamma^m\widehat\Gamma^{m\prime
}+
\widehat\Sigma_{\varepsilon\varepsilon}^m \bigr|+ m \frac{N\widebar K+T}{N\widebar
KT}\ln \bigl( \min(N
\widebar K,T) \bigr),
\]
where $\widehat\Gamma^m$ and $\widehat\Sigma_{\varepsilon\varepsilon}^m$ are
the respective estimators of $\Gamma$ and $\Sigma_{\varepsilon
\varepsilon}$ when the factor number is set to $m$ and $\widebar K=K+1$. In
the simulation, we set $r_{\max}=4$. For the model with zero
restrictions, we consider a two-step method to determine $r_1$~and~$r_2$. First, we use (\ref{eq0}) to estimate the total number
$r=r_1+r_2$, denoted by $\hat r$, and obtain $\hat\beta{}^{\hat r}$ by
the method of the basic model under $\hat r$. Then we calculate the
matrix $\mathscr{R}=(\mathscr{R}_{it})$ with $\mathscr R_{it}=\dot
y_{it}-\dot x_{it}\hat\beta{}^{\hat r}$ and use the information criterion
proposed by \cite{bai2002determining} to determine the factor number in
$\mathscr R$, which we use $\hat r_1$ to denote. In the second step,
the upper bound of the factor number is set to $\hat r$. Then $\hat
r_2=\hat r-\hat r_1$. For models in Section~\ref{sec4}, even though there are
observable common regressors and time invariant regressors in the $y$
equation, we treat them as part of the unknown factor structure when
estimating the total number of factors. Once the total number of
factors are obtained, the dimension of $g_t$ is obtained by subtracting
the dimension of $\phi_i$ and that of $d_t$ because $\phi_i$ and $d_t$
are observable in Section~\ref{sec4}. This approach works very well. Other
methods may also be considered.

We consider an unified way to estimate the model in Section~\ref{sec2} and the
model in Section~\ref{sec3} (with zero restrictions). More specifically, for a
given data set, we calculate $r$ and $r_1$. If $\hat r=\hat r_1$, we
turn to the basic model; if $\hat r>\hat r_1$, we turn to the model
with zero restrictions.

Tables~\ref{table1}--\ref{table2} report the simulation results based on 1000 repetitions.
Bias and root mean square error (RMSE) are computed to measure the
performance of the estimators. The percentage that the factor number is
correctly estimated by the above procedure is given in the third column
of each table. For comparison, we also report the performance of the
within-group (WG) estimators and Bai's iterated principal components
estimators (PC). Simulations for the models in Section~\ref{sec4} are provided
in the supplement~\cite{bai}.

\begin{sidewaystable}
\caption{The performance of WG, PC and ML estimators in the basic model}\label{table1}
\begin{tabular*}{\textwidth}{@{\extracolsep{\fill}}ld{3.0}d{3.1}ccccccccd{2.4}ccc@{}}
\hline
\multicolumn{1}{c}{$\bolds{N}$} & \multicolumn{1}{c}{$\bolds{T}$} & \multicolumn{1}{c}{\multirow{2}{20pt}{\centering{\textbf{\%} $\bolds{\hat r=r}$}}}
&\multicolumn{4}{c}{\textbf{WG}} & \multicolumn{4}{c}{\textbf{PC}} & \multicolumn{4}{c@{}}{\textbf{MLE}}
\\[-6pt]
& & &\multicolumn{4}{c}{\hrulefill} & \multicolumn{4}{c}{\hrulefill} & \multicolumn{4}{c@{}}{\hrulefill}
\\
& & & \multicolumn{2}{c}{$\bolds{\beta_1}$} & \multicolumn{2}{c}{$\bolds{\beta_2}$}
    & \multicolumn{2}{c}{$\bolds{\beta_1}$} & \multicolumn{2}{c}{$\bolds{\beta_2}$}
    & \multicolumn{2}{c}{$\bolds{\beta_1}$} & \multicolumn{2}{c@{}}{$\bolds{\beta_2}$}
\\[-6pt]
& & & \multicolumn{2}{c}{\hrulefill} & \multicolumn{2}{c}{\hrulefill}
    & \multicolumn{2}{c}{\hrulefill} & \multicolumn{2}{c}{\hrulefill}
    & \multicolumn{2}{c}{\hrulefill} & \multicolumn{2}{c@{}}{\hrulefill}
\\
& & & \multicolumn{1}{c}{\textbf{Bias}} & \multicolumn{1}{c}{\textbf{RMSE}} & \multicolumn{1}{c}{\textbf{Bias}}
    & \multicolumn{1}{c}{\textbf{RMSE}} & \multicolumn{1}{c}{\textbf{Bias}} & \multicolumn{1}{c}{\textbf{RMSE}}
    & \multicolumn{1}{c}{\textbf{Bias}} & \multicolumn{1}{c}{\textbf{RMSE}} & \multicolumn{1}{c}{\textbf{Bias}}
    & \multicolumn{1}{c}{\textbf{RMSE}} & \multicolumn{1}{c}{\textbf{Bias}} & \multicolumn{1}{c@{}}{\textbf{RMSE}}
\\
\hline
\phantom{0}50&75&99.9&0.1562&0.1616&0.1550&0.1600&0.0174&0.0405&0.0171&0.0411&-0.0001&0.0020&0.0000&0.0034\\
100&75&100.0&0.1539&0.1568&0.1558&0.1587&0.0061&0.0228&0.0062&0.0224&0.0000&0.0011&0.0000&0.0010\\
150&75&100.0&0.1534&0.1556&0.1540&0.1561&0.0029&0.0168&0.0028&0.0146&0.0000&0.0007&0.0000&0.0007
\\[3pt]
\phantom{0}50&125&100.0&0.1559&0.1605&0.1588&0.1636&0.0182&0.0389&0.0184&0.0409&0.0000&0.0017&0.0000&0.0016\\
100&125&100.0&0.1561&0.1586&0.1554&0.1579&0.0050&0.0167&0.0052&0.0167&0.0000&0.0009&0.0000&0.0008\\
150&125&100.0&0.1546&0.1565&0.1551&0.1570&0.0025&0.0108&0.0025&0.0106&0.0000&0.0006&0.0000&0.0005\\
\hline
\end{tabular*}\vspace*{20pt}
\caption{The performance of WG, PC and ML estimators in the model with zero restrictions}\label{table2}
\begin{tabular*}{\textwidth}{@{\extracolsep{\fill}}ld{3.0}d{3.1}cccccccccccc@{}}
\hline
\multicolumn{1}{c}{$\bolds{N}$} & \multicolumn{1}{c}{$\bolds{T}$} & \multicolumn{1}{c}{\multirow{2}{20pt}{\centering{\textbf{\%} $\bolds{\hat r=r}$}}}
&\multicolumn{4}{c}{\textbf{WG}} & \multicolumn{4}{c}{\textbf{PC}} & \multicolumn{4}{c@{}}{\textbf{MLE}}
\\[-6pt]
& & &\multicolumn{4}{c}{\hrulefill} & \multicolumn{4}{c}{\hrulefill} & \multicolumn{4}{c@{}}{\hrulefill}
\\
& & & \multicolumn{2}{c}{$\bolds{\beta_1}$} & \multicolumn{2}{c}{$\bolds{\beta_2}$}
    & \multicolumn{2}{c}{$\bolds{\beta_1}$} & \multicolumn{2}{c}{$\bolds{\beta_2}$}
    & \multicolumn{2}{c}{$\bolds{\beta_1}$} & \multicolumn{2}{c@{}}{$\bolds{\beta_2}$}
\\[-6pt]
& & & \multicolumn{2}{c}{\hrulefill} & \multicolumn{2}{c}{\hrulefill}
    & \multicolumn{2}{c}{\hrulefill} & \multicolumn{2}{c}{\hrulefill}
    & \multicolumn{2}{c}{\hrulefill} & \multicolumn{2}{c@{}}{\hrulefill}
\\
& & & \multicolumn{1}{c}{\textbf{Bias}} & \multicolumn{1}{c}{\textbf{RMSE}} & \multicolumn{1}{c}{\textbf{Bias}}
    & \multicolumn{1}{c}{\textbf{RMSE}} & \multicolumn{1}{c}{\textbf{Bias}} & \multicolumn{1}{c}{\textbf{RMSE}}
    & \multicolumn{1}{c}{\textbf{Bias}} & \multicolumn{1}{c}{\textbf{RMSE}} & \multicolumn{1}{c}{\textbf{Bias}}
    & \multicolumn{1}{c}{\textbf{RMSE}} & \multicolumn{1}{c}{\textbf{Bias}} & \multicolumn{1}{c@{}}{\textbf{RMSE}}
\\
\hline
\phantom{0}50&75&99.7&0.1098&0.1137&0.1095&0.1135&0.0097&0.0245&0.0099&0.0246&0.0000&0.0012&0.0000&0.0011\\
100&75&100.0&0.1088&0.1111&0.1092&0.1114&0.0038&0.0140&0.0038&0.0140&0.0000&0.0006&0.0000&0.0006\\
150&75&100.0&0.1086&0.1102&0.1083&0.1099&0.0011&0.0075&0.0015&0.0076&0.0000&0.0004&0.0000&0.0004
\\[3pt]
\phantom{0}50&125&99.7&0.1089&0.1121&0.1097&0.1130&0.0076&0.0199&0.0077&0.0196&0.0000&0.0009&0.0000&0.0009\\
100&125&100.0&0.1088&0.1107&0.1087&0.1106&0.0029&0.0104&0.0026&0.0100&0.0000&0.0005&0.0000&0.0004\\
150&125&100.0&0.1086&0.1099&0.1076&0.1090&0.0011&0.0055&0.0010&0.0054&0.0000&0.0003&0.0000&0.0003\\
\hline
\end{tabular*}
\end{sidewaystable}

From the tables, we can see that the factor number can be correctly
estimated with very high probability. It is also seen from the
simulations that the WG estimators are inconsistent. The bias of the WG
estimators shows no signs of decreasing as the sample size grows. The
iterated PC estimators are consistent, but biased. As the sample size
becomes large, the bias decreases noticeably. However, when the sample
size is moderate, the bias of the iterated PC estimators is still
pronounced. In comparison, the ML estimators are consistent and
unbiased. For all the sample sizes, the biases of the ML estimators are
very small and negligible. In addition, the RMSEs of the ML estimators
are always the smallest among the three estimators, illustrating the
efficiency of the ML method. The same patten is observed for all of the
four models considered.

\section{Conclusion}\label{sec7}
This paper considers estimating panel data models with interactive
effects, in which explanatory variables are correlated with the
unobserved effects. Standard panel data methods (such as the
within-group estimator) are not suitable for this type of models. We
study the maximum likelihood method and provide a rigorous analysis for
the asymptotic theory. While the analysis is difficult, the limiting
distributions of the MLE are simple and have intuitive interpretations.
The maximum likelihood method can incorporate parameter restrictions to
gain efficiency, a useful feature in view of the large number of
parameters under large $N$ and large $T$. We analyze the restrictions
via the Lagrange multiplier approach, which is capable of revealing
what kinds of restrictions lead to efficiency gain. We allow the model
to include time invariant regressors and common regressors. The
coefficients of the time invariant regressors are time dependent, and
the coefficients of the common regressors are cross-section dependent.
This is a sensible way for modeling the effects of such variables in
panel data context and fits naturally into the framework of interactive
effects. The likelihood method is easy to implement and performs very
well, as demonstrated by the Monte Carlo simulations.

\section*{Acknowledgments} The authors thank two anonymous referees, an Associate Editor and an
Editor for constructive comments.

\begin{supplement}[id=suppA]
\stitle{Supplement to ``Theory and methods of panel data models with interactive effects''}
\slink[doi]{10.1214/13-AOS1183SUPP} 
\sdatatype{.pdf}
\sfilename{aos1183\_supp.pdf}
\sdescription{This supplement provides detailed technical
proofs. Inferential theory for the estimated coefficients of
time-invariant and common regressors is given. The EM solutions are
shown to have local optimality property. Additional simulation results
are presented.}
\end{supplement}


\printaddresses

\end{document}